\documentclass{amsart}
\usepackage{amsfonts,amssymb,epsfig}
\usepackage[latin1]{inputenc}
\usepackage[francais]{babel}

\newtheorem{thm}{Théorème}[subsection]

\newtheorem{lemma}[thm]{Lemme}
\newtheorem{prop}[thm]{Proposition}
\theoremstyle{definition}
\newtheorem{defn}[thm]{Définition}
\theoremstyle{remark}

\theoremstyle{definition}
\newtheorem{example}[thm]{Example}

\newcommand{\bkR}{\mathbb R}

\begin{document}

\title{SUR LA CONJECTURE DES TORES PLATS}

\author{NGUIFFO B. BOYOM}
\address{GTA, UMR 5030 CNRS, Département de Mathématiques, Université Montpellier II}
\email{boyom@math.univ-montp2.fr}


\date{December 20, 2001}%

\begin{abstract} {The present work is devoted to compact completely solvable
solvmanifolds which admit Kahlerian metrics whose Kahler forms are homogeneous. In
particular, we show that such manifolds are diffeomorphic to flat tori. Our proof is
based on Dynkin diagrams associated to left invariant closed 2-forms in completely
sol\-va\-ble Lie groups. These diagrams are shown to be closely related to
Lagrangian fo\-lia\-tions,primitive degrees of homogeneous spaces , and geometric
structures such as $ (\bkR^ n, Aff (\bkR^ n ) ) $-structures on symplectic
homogeneous spaces. \par Our result may be viewed as a non trivial generalization of
a result by C. Benson ans C. Gordon, $\lbrack$BG1$\rbrack$, who showed that compact
nilmanifolds $ \Gamma \backslash G $ that admit Kahlerian metrics are diffeomorphic
to flat tori. Their proof is based on rational homotopy theory of compact Kahlerian
manifolds $\lbrack$DGMS$\rbrack$, the hard Lefschetz theorem and a theorem of K.
Nomizu which says that $\bigoplus_k H^k_{De\ Rham}(\Gamma \backslash G) $ is
isomorphic to the algebra $\bigoplus_k H^k(L(G),\bkR) , L(G) $ being the Lie algebra
of G. McDuff $\lbrack$Mc D1$\rbrack$ gave another proof of Benson-Gordon theorem, by
using moments maps of symplectic $ S^1 $-actions. The main starting arguments in the
both proofs have no equivalent in the case of general solvable Lie group (e.g. the
center of a nilpotent Lie group contains a nontrivial 1-parameter subgroup). To work
with completely solvable Lie group, we use very different approach from the ones
used by Benson-Gordon and McDuff, (e.g. the existence of affinely flat bilagrangian
structures in completely solvable symplectic homogeneous spaces
$\lbrack$NB4$\rbrack$, $\lbrack$NB5$\rbrack$, Marsden-Weinstein reduction theory,
$\lbrack$MW$\rbrack$ }
\end{abstract}
\maketitle

\section{Introduction}

Dans ce travail un groupe discret  $\Gamma$ est appelé groupe Kahlérien s'il est
isomorphe au groupe fondamental d'une variété analytique complexe admettant une
métrique Kählérienne . Dans ce travail nous présentons une variété Kählérienne sous
la forme de triplet $( M, \omega , J )$ où $\omega$ est la forme de Kähler et J est
le tenseur de la structure presque complexe de la variété analytique réelle M . On
ne sait pas caractériser les groupes fondamentaux des variétés Kählériennes
compactes .Il existe cependant des situations où l'exsitence de type particulier de
système dynamique différentiable sur M determine la géométrie globale de M , c'est à
dire son type de difféomorphisme . Il en est ainsi quand M est une nilvariété
compacte . La première démonstration publiée du théorème suivant est due à C. BENSON
et C. GORDON dans un travail commun .

\begin{thm}[\cite{BG 1}] Soit  $G $ un groupe de Lie nilpotent connexe , simplement connexe
, contenant un réseau $\Gamma$ . Si la variété analytique réelle $\Gamma \backslash
G $ admet une métrique Kählérienne alors $\Gamma \backslash G $ est difféomorphe au
tore plat.
\end{thm}

De ce théorème Dusa McDUFF a publié une seconde démonstration en 1978 . La
conjecture des tores plates concerne l'analogue du théorème ci-dessus pour les
groupes de Lie complètement résolubles , un des objectifs de ce travail est d'en
donner la démonstration . Si $ G $ est un groupe de Lie connexe , simplemet connexe
et contenant un réseau cocompact $\Gamma$ alors une k-forme différentielle définie
sur la variété quotient $\Gamma\backslash G $ sera dite G-homogène (ou tout
simplement homogène quand il n'y a pas de risque de confusion ) si son image inverse
par la projection canonique de $ G $ sur $\Gamma\backslash G $ est une k-forme
invariante par les translations à gauche définies par les éléments de $ G $ . Nous
démontrerons à la section 4 de ce travail le résultat suivant.

\begin{thm} Soit $ G $ un groupe de Lie complètement résoluble connexe et
simplement connxe . On suppose que $ G $ contient un réseau cocompact $\Gamma$ et
que la variété analytique réelle $\Gamma\backslash G $ admet un métrique Kählérienne
dont la forme de Kähler est homogène . Alors la variété analytique réelle
$\Gamma\backslash G $ est difféomorphe au tore plat .
\end{thm}

Du point de vue de la topologie de la variété réelle $\Gamma\backslash G $
l'hypothèse de l'homogenéité de la forme de Kähler n'est pas vraiment une
restriction . En effet la complètement résolubilité du groupe de Lie $ G $ entraine
que l'algèbre de De RHAM de la variété réelle $\Gamma\backslash G $ est isomporphe à
l'algèbre de KOSZUL de l'algèbre de Lie du groupe de Lie $ G $ ; cela est une
généralisation d'un théorème classique de K Nomizu pour les nilvariétés compactes
[NK] . Cette généralisation est le fruit de suite spectrale consultable soit dans A.
Hattori [HA] soit dans A. RAGHUNATHAN [RA]. Il est opportun de signaler que la
démonstrations de BENSON-GORDON  n'est pas adaptables au cas des groupes de Lie
complètement résolubles, bien que les démonstrations utilisent le théorème de
Nomizu. En effet cette démonstration s'appuie dès le départ sur des propriétés des
groupes de Lie nilpotents qui n'ont pas d'équivalents dans les groupes de Lie
complètement résolubles ,elle fait intervenir des propriétés topologiques profondes
des variétés Kählériennes compactes , e.g. homotopie rationnelle [DGMS], théorème
dur de LEFSCHEZ, triple produit de MASSEY. La même réserve vaut pour la
démonstration de Dusa McDUFF qui combine une version affaiblie du théorème dur de
LEFSCHEZ avec l'application moment généralisé des actions hamiltoniennes des cercles
( le cercle utilisé vient d'un sous-groupe à un paramètre central intersectant le
réseau $\Gamma$ suivant un réseau propre ) . Il convient de signaler que les
techniques utilisées aussi bien dans [BGl] que dans [McDl] ont permis la découverte
de nombreux exemples de variété symplectiques compactes n'admettant pas de métrique
de Kähler . Il y une literature abondance dans cette direction, [CFG], [McD2], [OT]
et c..  Dans un travail à part [ NB6 ] nous avons mis au point des concepts qui
facilitent ici des incursions efficaces de la géométrie affine (voir tout
particulièrement les sections 2 et 3 de ce travail ).

Les critiques et les remarques du reféré ont permis d'améliorer la présentation de
concepts introduits dans [ NB6 ] ainsi que leur utilisation dans ce travail .
L'auteur lui en sait gré .

\section{Preliminaires}

\subsection{Diagrammes des 2-formes différentielles}

Soit $G$ un groupes de Lie connexe , simplement connexe et complètement résoluble.
L'algèbre de Lie des champs de vecteurs invariants par les translations à gauche sur
$G$ est notée $\frak g$ . On désigne par   $\Omega_l^2(G)$ l'espace vectoriel des
2-formes différentielles fermée et invariantes par les translations à gauche . Nous
allons rappeler brièvement la notion de diagramme attaché à un couple  $(\omega , F
) \in \Omega^2_l( G ) \times {\mathcal F}(G )$ où   ${\mathcal F}( G )$ est
l'ensemble des séries de composition du groupe de Lie $G$ . Une série de composition
est une filtration de $G$ par des sous-groupes de Lie $G_k$ tels que $\dim\ G_k = k$
et $G_k$ est distingué dans $G_{k+1}$. Soit $i_k$ l'homomorphisme inclusion de $G_k$
dans $G$ et $\omega_k = i_k^{\ast}\omega$. Soit $g_k$ l'algèbre de Lie de $G_k$ et
$h_k$ le noyau de la restriction à $g_k$ de $k$ . On désigne par $H_k$ le
sous-groupe de Lie connexe de $G_k$ déterminé par $h_k$ , alors les sous-groupes de
Lie $H_k$ et $H_{k+1}$ sont reliés par une des relations d'inclusion soit $H_k
\subset H_{k+1}$ soit $H_{k+1} \subset H_k$ . Lorsque $H_k \subset H_{k+1}$ les
variétés symplectiques $( G_k/H_k , \omega_k )$ et $( G_{k+1}/H_{k+1} , \omega_{k+1}
)$ sont canoniquement isomorphes. Lorsque $H_{k+1}$ est inclu dans $H_k$ la variété
symplectique $( G_k/H_k , \omega_k )$ est une réduite d'une action hamiltonienne du
groupe $\bkR$ dans $(G_{k+1}/H_{k+1} , \omega_{k+1} )$, [ MW ] . A tout couple
$(\omega , F ) \in \Omega_l^2( G ) \times {\mathcal F} ( G )$ on fait correspondre
un diagramme $d( \omega, F )$ dont les sommets ont pour coordonnées les couples $(
G_k , H_k )$ . Ce diagramme est obtenu par adjonction bout à bout des morceaux
semblables aux modèles élémentaires suivants

$$
\circ \rightarrow \circ \leftrightarrows \circ \ \ , \ \ \circ \leftrightarrows
\circ \rightarrow \circ
$$

On a introduit dans  [NB6] les notions suivantes : diagramme semi-normal ; diagramme
semi-nilpotent ; diagramme simple . Pour le travail présent on retiendra que les
deux affirmations suivantes sont équivalentes:

(a) la forme  $\omega$ possède un diagramme simple   $d (  \omega , F )$

(b) la variété symplectique $( G/H ,  \omega )$ possède un feuilletage lagrangien L
qui est invariant par l'action à gauche de $G$ sur $G/H$ ; $H$ est le sous-groupe de
Lie connexe déterminé par le noyau de  $ \omega : \frak g \times \frak g \rightarrow
R$ . On a démontré dans [NB6] que si une 2-forme   $\omega \in \Omega_l^2 ( G )$
possède un diagramme qui est semi-nilpotent et semi-normal alors elle possède un
diagramme simple . Cette notion fournit un cadre pour traiter les formes de Kähler
qui apparaissent dans la conjecture des tores plats . Ce fait a été mis en évidence
par C. Benson et C. Gordon .

\begin{thm} [ \cite{BG 2} ] . Soit $G$ un groupe de Lie connexe , simplement connexe et complètement
résoluble . On suppose que $G$ contient un réseau cocompact  $\Gamma$ et que $\omega
\in \Omega_l^2 ( G )$ est une forme de rang maximal . Si la variété symplectique $(
G , \omega)$ possède une métrique Kählérienne $h = g +  \sqrt{-1}\omega$    alors le
sous-groupe des commutateurs $[G,G]$ hérite de $\omega$  d'une forme symplectique .
\end{thm}

\subsection{Formes différentielles cocompactes}

 Dans cette sous-section il sera question des variétés homogènes symplectiques
compactes dont les formes symplectiques sont homogènes.  Soit $ G $ un groupe de Lie
connexe et simplement connexe quelconque. Soit $ \omega  \in  \Omega^ 2_{l}( G) $ ;
le sous-groupe de Lie connexe $ H_\omega $ dont l'algèbre de Lie engendre le noyau
de $ \omega $ est supposé fermé dans $G$.

\begin{defn}
Un sous-groupe (abstrait) $ \Gamma \subset G $ vérifie la condition $ (c_\omega ) $
si l'espace des orbites $ \Gamma \backslash G/H_\omega $ de l'action à gauche de $
\Gamma $  dans $ G/H_\omega $ est une variété compacte de dimension positive.
\end{defn}

Dans la suite, on désignera par $ {\mathcal C}_\omega( G) $ la collection des
sous-groupes $ \Gamma $ de $ G $ qui vérifie la condition $ (c_\omega )\, . $

\begin{defn} Une 2-forme  $\omega \in \Omega^ 2_{\ell}( G) $  est dite cocompacte si
$ {\mathcal C}_\omega( G) $  contient un sous-groupe discret $ \Gamma_ \omega . $
\end{defn}

Soit $ \omega  \in  \Omega^ 2_{\ell}( G) $ une forme cocompacte ; les deux variétés
homogènes $ H_\omega \backslash G $ et $ G/H_\omega $ possèdent des structures
symplectiques héritées de la forme $ \omega \, . $ Visiblement $ (G/H_\omega ,\omega
) $ est $ G $-homogène. Quand une variété symplectique $ (M,\omega) $ possède un
feuilletage lagrangien on dira pour abréger qu'elle est $ L ${\sl -feuilletée.\/}
Puisque $ \omega $ est supposé \^etre cocompacte, fixons un sous-groupe discret $
\Gamma_ \omega \in  {\mathcal C}_\omega( G)\, . $ Le difféomorphisme inversion $
{\rm Inv}(g) = g^{-1} $ se factorise pour définir un difféomorphisme $ \varphi $ de
$ \Gamma_ \omega \backslash G/H_\omega $ sur $ H_\omega \backslash G/\Gamma_ \omega
. $ Par ailleurs la forme symplectique de $ G/H_\omega $ héritée de $ \omega$  se
projette en une forme symplectique dans $ \Gamma_ \omega \backslash G/H_\omega ; $ $
H_\omega $ étant connexe $ G/H_\omega $ est le revêtement universel de $ \Gamma_
\omega \backslash G/H_\omega , $ l'application canonique de $ G/H_\omega $ sur $
\Gamma_ \omega \backslash G/H_\omega $ est un revêtement symplectique

Supposons $ G $ complètement résoluble, le théorème 1.4.1 de [NB6] assure que si $
\omega $ possède un diagramme pondéré simple et est cocompacte, alors $ G/H_\omega $
et $ \Gamma_ \omega \backslash G/H_\omega $ sont $ L $-feuilletées. Le revêtement
universel de $ G/H_\omega $ de $ \Gamma_ \omega \backslash G/H_\omega $ possède une
structure affine dont l'holonomie linéaire est inclue dans $ S\ell( m,\bkR) \, , $ $
m = \dim \Gamma_ \omega \backslash G/H_\omega \, . $

Supposons que $ \omega  \in  \Omega^ 2_{\ell}( G) $ soit symplectique, alors elle
est cocompacte si et seulement si $ \Gamma_ \omega $ est un réseau cocompact dans $
G\ $

\noindent On observe que si $ \Gamma $ est un réseau cocomppact dans $ G\, et si , $
$ \Gamma \backslash G $ possède une 2-forme symplectique alors elle en possède une
autre forme symplectique qui est $ G $-homogène au sens de [McD1] , c'est-à-dire qui
se relève en une forme symplectique invariante à gauche dans $G$, [HA], mais $
(\Gamma \backslash G,\omega) $ n'est pas en général $ G $-homogène.

\subsection{Déformations des structures kähleriennes
$(M,\omega,J)$ qui sont  $\Gamma_{\omega}$-invariantes.}

Nous revenons au cas des groupes de Lie complètement résolubles. Soit $ G $ un tel
groupe. Si $ \omega \in \Omega^ 2_{\ell}( G) $ est cocompacte suivant la définition
1.2.2 alors on fixe un sous-groupe $\Gamma_\omega$   dans $C_\omega  (G )$; $ M $
désigne la variété $ G/H_{\omega} $ munie de l'action à gauche de $ \Gamma_ \omega
\in {\mathcal C}_\omega( G),\Gamma_ \omega $ est discret.

Soit $ {\mathcal K}^a_\omega( M) $ le ``module sur $\Gamma_ \omega$'' des structures
k\"ahlériennes affinement plates $ (M,\omega ,J,D) $ dans $ (M,\omega) \, ; $ on
désigne par $ (M,\omega ,J,D) $ la structure Kählérienne $ (M,\omega ,J) $ reliée à
la structure affinement plate $ (M,D) $ par les deux relations $ D\omega  = 0 $ et $
DJ = 0\, . $ \par \noindent En vertu des résultats mis en évidence dans [NB6] on
sait que $ {\mathcal K}^a_\omega( M) $ est non vide quand la combinatoire de $
\omega $ est raisonable. Il en est ainsi en particulier quand $ \omega $ possède un
diagramme pondéré (ou un graphe orienté) simple $ {\rm O} \rightarrow {\rm
O}^{\leftarrow}_{ \rightarrow} {\rm O} \, . $ \par
\bigskip
\noindent{\bf Question 1.} {\sl Existe-t-il une structure \/}$ (M,\omega ,J,D)\in
{\mathcal K}^a_\omega( M) $ {\sl qui est invariante par \/}$ \Gamma_ \omega ? $ \par
\smallskip
\noindent Rappelons que la conjecture des tores plats tente de répondre à la
question 2 suivante : $ \omega  \in  \Omega^ 2_{\ell}( G) $ est symplectique et
cocompacte. \par
\bigskip
\noindent{\bf Question 2.} {\sl On suppose que la variété symplectique \/}$ (\Gamma_
\omega \backslash G,\omega ) $ {\sl possède une structure k\"ahlérienne \/}$
(\Gamma_ \omega \backslash G,\omega ,J )\, , $ {\sl à quoi ressemble la variété
compacte \/}$ \Gamma_ \omega \backslash G\, . $

Naturellement si on note $ \pi $ la projection de $ M = G/H_\omega $ sur $ \Gamma_
\omega \backslash G/H_\omega $ et $ {\mathcal K}_\omega $ l'ensemble des structures
k\"ahlériennes dans $ (\Gamma_ \omega \backslash G/H_\omega )\, , $ toute structure
$ (G_\omega \backslash G/H_\omega ,\omega ,J ) $ se relève en une structure
Kählérienne  $(M,\omega ,J) $ dans $ (M,\omega)$; notons $ {\mathcal K}_\omega( M) $
les structures ainsi obtenues dans $(M,\omega)$ (qui se projettent dans $ \Gamma_
\omega \backslash G/H_\omega$). La question 1 équivaut à l'examen de $ {\mathcal
K}^a_\omega( M) \cap {\mathcal K}_\omega( M)$.

Dans un premier temps nous allons démontrer le théorème clef suivant

\begin{thm}
\label{thm:231} Si $ \omega $  possède un diagramme simple $ d(\omega ,F) $ alors si
$ {\mathcal K}_\omega( M) $  est non vide alors $ {\mathcal K}_\omega( M) \cap
{\mathcal K}^a_\omega( M) $  est non vide aussi
\end{thm}

La démonstration du théorème \ref{thm:231} s'appuie sur le théorème 1.4.1 de [NB6]
joint au théorème de BRUHAT--WHITNEY\footnote{voir Théorème \ref{thm:324}} suivant

\begin{thm}[\cite{BW}]. Soient $ (M,J) $ et $ (M,J^\prime) $ deux structures de variétés
analytiques complexes ayant la même structure de variété analytique réelle
sous-jacente. Soit $ N $  une sous-variété analytique réelle de $ M $  de dimension
$ {1 \over 2} \dim M\, . $
\end{thm}

Si $ N $  est totalement réelle pour les deux structures complexes $ (M,J) $  et
\par \noindent$ (M,J^\prime) \, , $ alors il existe un voisinage ouvert $ U $  de $
N $ dans $ (M,J)\, , $ un voisinage ouvert $ U^\prime $ de $ N $  dans $
(M,J^\prime) $  et un difféomorphisme analytique $ \varphi $  de $ U $ sur $
U^\prime $ satisfaisant les conditions (i) $ \varphi( x) = x $  pour tout $ x \in  N
$ et (ii) $ d\varphi_ 0J = J^{\prime}_ od\varphi \ $ \\

{\it Démonstration du théorème \ref{thm:231}}. Nous supposons que la 2-forme $
\omega \in \Omega^ 2_{\ell}( M) $ est cocompacte. Dans toute la suite on pose $ H =
H\omega \, . $ La variété compacte $ \Gamma_ \omega \backslash G/H $ a pour
revêtement universel la variété $ M = G/H\, . $ De plus ce revêtememnt est
symplecltique. La démonstration se fait en cinq étapes. On pose $ 2m = \dim M\, . $
\par

\noindent{\bf 1ère étape.} Puisque $ \omega $ possède un diagramme pondéré simple le
théorème 1.4.1 de [NB6 ] assure l'existence dans la variété symplectique $
(M,\omega) $ d'une structure bilagran\-gienne $ (L,{\mathcal N}) $ dont le
feuilletage $ L $ est $ G $-invariant. De plus la paire $ (L,{\mathcal N}) $ donne
lieu à une unique structure affinement plate $ (M,D) $ qui se comporte vis à vis de
$ (M,\omega ,L,{\mathcal N}) $ suivant les conditions ci-dessous :

\smallskip
\noindent (i) $ DL \subset  L $ et $ D{\mathcal N}\subset{\mathcal N}\, . $ \par
\smallskip
\noindent (ii) $ D\omega  = 0\, . $ \par
\smallskip
\noindent (iii) Dans un voisinage de chaque point $ x\in M\, , $ il existe un
système de coordonnées de Darboux $ (q_1,...,q_m,p_1,...,p_m ) $ compatible avec $
(L,{\mathcal N}) $ dans le sens suivant : \par \noindent les fonctions locales $ q_i
$ sont des intégrales premières de $ {\mathcal N} $ et les fonctions locales $ p_i $
sont des intégrales premières de $ L\, . $ \par
\smallskip
\noindent (iv) Soit $ (q^{\prime}_ 1,...,q^{\prime}_ m,p^{\prime}_ 1,...,p^{\prime}_
m ) $ un autre système des coordonnées de Darboux définies dans un voisinage du même
point $ x\in M\, , $  et satisfaisant la condition (iii ) alors ces coordonnées
s'expriment affinement en fonction des précédentes $$ q^{\prime}_ i = \sum^
m_{j=1}A_{ij}q_j+ c_i, p^{\prime}_ i = \sum^ m_{j=1}A_{ij}p_j+ d_i\, , $$ La matrice
$ A = [A_{ij} ]\, , $ $ 1 \leq  i,j \leq  m $ étant orthogonale.
\par
\medskip
\noindent Remarques. (a) On vérifie sans difficulté que tout système des coordonnées
de Darboux qui vérifie (iii) est aussi un système des coordonnées affines de la
structure localement plate $ (M,D)\, , $ (voir $\lbrack$NB j], $ j = 4,5). $
\par
\smallskip
\noindent (b) La condition (iv) entra\^\i ne que $ M $ porte une structure de
variété analytique complexe dont des coordonnées locales sont les

\begin{equation}
z_i= q_i+ \sqrt{ -1} p_i, 1 \leq  i \leq  m\, .
\end{equation}

(c) La même condition (iv) assure que munie de l'atlas complexe défini par les
coordonnées $ z_i $ ci-dessus, la variété complexe $ M $ porte une structure
hermitienne $ h $ qui s'exprime dans les coordonnées $ z_i $ par

\begin{equation} h = \sum^{ }_ i dz_id\bar  z_i\, .
\end{equation}

(d) visiblement $ h $ est plate et la remarque (a) ci-dessus montre que $ D $ est la
connexion de Levi--Civita de $ h\, . $ \par
\bigskip
\noindent{\bf Conclusion partielle 1.} La variété symplectique $ (M,\omega) $ porte
une structure k\"ahlérienne affinement plate $ (M,\omega ,J)\, . $ Les feuilles des
feuilletages $ L $ et $ {\mathcal N} $ sont des sous-variétés totalement géodésiques
de $ (M,D) $ et totalement réelles de $ (M,J)\, $ \par
\bigskip
\noindent{\bf 2ème étape.} Dorénavant un système de coordonnées locales qui
satisfait les conditions (iii) et (iv ) de  (1) est appelé coordonnées de
Darboux-Hess. On a supposé que $ {\mathcal K}_\omega( M) $ est non vide ; ce qui
équivaut à dire que la variété symplectique compacte $ \Gamma_ \omega \backslash M =
\Gamma_ \omega \backslash G/H $ porte une structure k\"ahlérienne $ (\Gamma_ \omega
\backslash M,\omega ,J_0 )\, . $
\par
\noindent On note encore $ (M,\omega ,J_0 ) $ la structure k\"ahlérienne image
inverse de $ (\Gamma_ \omega \backslash G/H,\omega ,J_0 ) $ par l'application
revêtement $ M \mapsto  \Gamma_ \omega \backslash M\, . $ Les feuilles des
feuilletages lagrangiens $ L $ et $ {\mathcal N} $ sont de nouveau des sous-variétés
totalement réelles de $ (M,J_0 )\, . $
\par \noindent On observe que des deux feuilletages $ L $ et $ {\mathcal N} $ le premier
$ L $ est invariant par $ \Gamma_ \omega \, , $ qui est le groupe fondamental de $
\Gamma_ \omega \backslash M\, . $ \par
\bigskip
\noindent{\bf 3ème étape}\footnote{Voir Théorème \ref{thm:324}} Soit $ L_x $ la
feuille de $ L $ qui passe par le point $ x \in  M\, . $ En vertu du théorème de
BRUHAT--WHITNEY, on va choisir deux voisinages ouverts $ U $ et $ U^\prime $ de la
sous-varité $ L_x $ et un difféomorphisme $ \varphi $ de $ U $ sur $ U^\prime $ qui
prolonge le difféomorphisme identité de $ L_x $ et satisfait la condition
$$ d\varphi  \circ  J = J_0\circ  d\varphi $$
Soit $ (q_1,...,q_m,p_1,...,p_m ) $ un système des coordonnées de Darboux définies
dans un voisinage $ x $ contenu dans $ U \cap  U^\prime $ et vérifiant les
conditions (i) à (iv) de la première étape. On suppose que la feuille $ L_x $ de $ L
$ est définie par $ p_j= 0\, , $ $ 1 \leq  j \leq  m\, . $ Dans les coordonnées $
(q_i,p_i ) $ comme ci-dessus. \par \noindent On pose $ (q,p) =
(q_1,...,q_m,p_1,...,p_m )\, ; $ alors on a un difféomorphisme holomorphe local dont
une forme normale est la suivante
$$ \varphi( q,p) = (q+h(p)\alpha( q,p),k(p)\beta( q,p) $$
o\`u les fonctions $ h $ et $ k $ vérifient
$$ h(0) = k(0) = 0\, . $$
La matrice de $ d\varphi $ dans les coordonnées $ (q,p) $ a la forme suivante

$$ d\varphi( q,p) = [ \begin{pmatrix}1 + h(p)a(q,p) &
b(q,p) \\ k(p)c(q,p) &  d(q,p) \end{pmatrix} ]
$$

où $ h(p) $ et $ k(p) $ sont des fonctions dont les 1-jets sont nuls pour $ p = 0\,
. $ \par \noindent Notons $ TL $ et $ T {\mathcal N} $ les distributions tangentes à
$ L $ et à $ {\mathcal N} $ respectivement. On a bien entendu

$$ JTL = T{\mathcal N}\, , $$

par conséquent la relation

\begin{equation}
d\varphi  \circ  J = J_o\circ  d\varphi
\end{equation}

montre que $ d\varphi_ x $ envoie $ T_x{\mathcal N} $ sur le sous-espace orthogonal
du sous-espace $ T_xL $ relativement à la métrique riemannienne de $ (M,\omega ,J_0
)\, . $ Soient $ X $ et $ Y $ deux sections locales de $ TL\, , $ alors $ JX $ et $
JY $ sont des sections locales de $ T{\mathcal N} $ qui sont envoyés par $ d\varphi_
y $ dans $ J_0T_yL $ pour tout $ y \in L_x\, . $ La relation

\begin{equation}
 d\varphi_ y[JX,JY]_y= [d\varphi JX,d\varphi JY]_y, y \in  L_x
\end{equation}

montre que le système différentiel $ J_0TL $ est intégrable le long de $ L_x\, . $
En d'autres termes par tout $ y \in  L_x $ passe une germe de sous-variété intégrale
maximale de $ J_0TL $ de dimension $ m $ qui n'est pas autre chose que le germe en $
y $ de l'image par $\phi $ de la feuille $ {\mathcal N}_y $ de $ {\mathcal N}\, . $

Nous concluons que le système différentiel $ J_0TL $ est complètement
intégrable\footnote{Voir aussi Théorème \ref{thm:324}}
\par

\noindent{\bf Conclusion partielle 2.} On a déduit de la structure bilagrangienne $
(L,{\mathcal N}) $ la structure bilagrangienne $ (L,J_0L )\, . $ Puisque $ L,J_0 $
et $ \omega $ sont invariants par $ \Gamma_ \omega $ (ou si on veut par le groupe
fondamental de $ \Gamma_ \omega \backslash M), $ $ (L,J_0L ) $ se projette en une
structure bilagrangienne dans la variété symplectlique $ (\Gamma_ \omega \backslash
M )\, . $

Pour simplifier on pose $ J_0L $ pour $ J_0TL\, . $ \\

\noindent{\bf 4ème étape.} On note $ D° $ la connexion linéaire sans torsion définie
dans $ M $ par la paire de feuilletages lalgrangiens $ (L,J_0L )\, , $ (voir (17)).
On a $ D°L \subset  L\, , $ $ D°J_0L \subset  J_0L $ et $ D°\omega  = 0\, . $
\par
\bigskip
Remarque . La connexion linéaire $D°$ et la connexion $D$ définie  par la paire
$(L,N )$ coincident sur $L$. Ce sont des connexions symplectiques préservant le
feuilletage lagrangien $L$ ,étant toutes les deux de torsion nulle il diffèrent par
un tenseur symétrique dont la valeur en chaque point est dans le premier
prolongement à la Guillemin-Sternberg de $Sp(\omega)$ ,  ( [ GS ] , [SS] )

\begin{lemma}
\label{lem:232} Le tenseur de courbure de $ D° $  est nul.\footnote{Voir aussi
Théorèmes \ref{thm:326} et \ref{thm:327}}
\end{lemma}

\noindent Preuve. On considère un système des coordonnées de Darboux $ (q_i,p_i ) $
de $ (M,\omega) $ satisfaisant les conditions (i) à (iv) de la première étape.
\par \noindent On pose $ X_i= X_{q_i} $ et $ Y_i= X_{p_i} $ les champs hamiltoniens
correspondant. On a $ JX_i= Y_i $$\quad$ $ 1 \leq  i \leq  m\, , $ $ 2m = \dim M\, .
$ Soit $ \varphi $ le difféomorphisme local du théorème de BRUHAT--WHITNEY. On a

\begin{equation}
J_0X_i(y) = d\varphi_ yY_i(y)
\end{equation}

pour tout $ y \in  L_x\, . $

Le tenseur de courbure de $ D $ est nul sur $ L \times L $ et sur $ J_0L \times
J_0L\, . $ Il reste à calculer ce tenseur pour une par $ (X_i,J_0X_j )\, . $ \par
\noindent Auparavant on observe que $ D $ co\"\i ncide avec $ D^0 $ aussi bien sur $
L \times  L $ que sur $ J_0L \times J_0L\, . $ Autrement dit
$$ D_{X_i}X_j= D^0_{X_i}X_j\, . $$
Et puisque les $ X_i $ sont hamiltoniens et engendrent $ L\, , $ on a

\begin{equation}
 i (Di^0_{X_i}X_j )\omega  = L _{X_i}i (X_j )\omega  = 0\, ,
\end{equation}

c'est à dire que $ D_{X_i}X_j= D^0_{X_i}X_j= 0=. $ Puisque le tenseur de courbure $
R $ de $ D $ vaut
$$ R (X_i,J_oX_j ) = D_{X_i}D_{J_oX_j}- D_{J_oX_j}D_{X_i}- D_{
[X_i,J_oX_j ]}\, , $$ et compte tenu de $ D_{X_i}X_j = 0 $ on a
$$ R (X_i,J_oX_j )X_k= D_{X_i}D_{J_oX_j}X_k - D_{ [X_iJ_oX_j
]}X_k\, . $$ Les calculs ci-dessous sont effectués en un point $ y $ de $ L_x $
qu'on omettra volontairement d'expliciter. Ainsi au pont $ y \in  L_x $ on a

\begin{equation}
\omega (R (X_i,J_oX_j )X_k,Y_{\ell} ) = \omega (D_{X_i}D_{jX_J}X_k- D_{ [X_iJ_oX_j
]}X_k,Y_{\ell} )
\end{equation}

On tient maintenant compte du fait que $ DL \subset L$; d'où

$$
\begin{array}{rl}
D_{X_i}D_{J_oX_j}X_k & = D^o_XD_{J_oX_j}X_k\ {\rm avec} \\
i (D^o_XX^\prime )\omega &  = L(X)L(X^\prime) \omega \ {\rm pour} \ X\ {\rm et} \
X^\prime \ {\rm tangent} \ {\rm \grave a} \ L\, .
\end{array}
$$

L'expression (20) devient maintenant

$$
\begin{array}{ll}
\omega (R (X_i,J_oX_j )X_k,Y_{\ell} ) = & X_i\omega (D_{J_oX_j}X_k,Y_{\ell} ) -
\omega (D_{J_oX_j}X_k, [X_i,Y_{\ell} ] ) \\  &  - \omega ( [ [X_i,J_oX_j ],X_k ]-
D_{X_k} [X_i,J_oX_j ],Y_{\ell} )
\end{array}
$$

Puisque $ [X_i,Y_{\ell} ] = 0\, , $ l'expression ci-dessus est réduite à
$$ X_i\omega (D_{J_oX_j}X_k,Y_{\ell} ) + \omega ( [X_k,
[X_i,J_oX_j ] ],Y_{\ell} ) - \omega (D_{X_k} [X_i,J_oX_j ],Y_{\ell} )\, . $$ En
vertu de $ D_{X_k}\omega  = 0\, , $ la dernière expression devient

$$
\begin{array}{ll}
\omega (R (X_i,J_oX_j )X_k,Y_{\ell} ) & = X_i\omega (D_{J_oX_j}X_k,Y_{\ell} ) + X_k
(X_i\omega (J_oX_j,Y_{\ell} ) ) \cr  &  - X_k (X_i\omega (J_oX_j,Y_{\ell} ) ) +
\omega ( [X_i,J_oX_j ],D_{X_k}Y_{\ell} ) \\  &  = X_i\omega (D_{J_oX_j}X_k,Y_{\ell}
) + \omega ( [X_i,J_oX_j ],D_{X_k}Y_{\ell} )\, .
\end{array}
$$

Etant donné que $ L (X_i )\omega  = 0\, , $ on voit que

$$
\begin{array}{l}
X_i\omega (D_{J_oX_j}X_k,Y_l) + \omega ([X_i,J_oX_j],D_{X_k}Y_l) \\
= X_i\omega (D_{J_oX_j}X_k,Y_l) + X_i\omega (J_oX_j,D_{X_k}Y_l) - \omega
(J_oX_j,[X_i,D_{X_k}Y_l]) \\
= X_i\omega ([J_oX_j,X_k],Y_l) + X_i\omega (D_{X_k}J_oX_j,Y_l) + X_i\omega
(J_oX_j,D_{X_k}Y_l) \\
- \omega (J_oX_j,[X_i,D_{X_k}Y_l]) \ .
\end{array}
$$

Les arguments précédents i.e. $ L_{X_i}\omega  = 0, \, D\omega  = 0, \,
[X_i,Y_{\ell} ] = 0 $ montrent que l'expression ci-dessus se réduit à

$$
\begin{array}{rl}
\omega (R (X_i,JX_j )X_k,Y_{\ell} ) = & - X_i (X_k\omega (J_oX_j,Y_{\ell} ) + X_i
(X_k\omega (J_oX_j,Y_{\ell} ) ) )
\\  &  - \omega (J_oX_j, [X_i,D_{X_k}Y_{\ell} ] )
\end{array}
$$

Puisque le tenseur de torsion de $ D $ est nul et que $ [X_k,Y_{\ell} ] = 0 $ on a

$$ \begin{array}{rl} \omega (R (X_i,JX_j )X_k,Y_l ) & = -
\omega (J_oX_j, [X_i,D_{Y_l} X_k ] ) \\  &  = - \omega (J_oX_j,D_{X_i}D_{Y_l} X_k- D
_{D_{Y_l} X_k}X_i ) \\  &  = \omega (J_oX_j,D_{D_{Y_l} X_k}X_i- D_{Y\ell} D_{X_i}X_k
) \\  &  = w (J_oX_j,D^o_{D_{Y_l} X_k}X_i )
\end{array}
$$

Puisque $ D^o_XX_i= 0 $ pour tout $ X $ tangent à $ L $ on conclut que

$$ \omega (R (X_i,J_oX_j )X_k,Y_l ) = 0 $$

d'où

$$ R (X_i,J_oX_j ) = 0\quad 1 \leq  i,j \leq  m\, . $$

Cela termine la démonstration du Lemme \ref{lem:232}. \\

\noindent{\bf Conclusion partielle 3.} La paire $ (L,J_oL ) $ détermine une
structure bilagran\-gienne affinement plate dans $ (M,\omega) $ qui est invariante
par l'action de $ \Gamma_ \omega \, . $ Il donne donc naissance à une structure
affine plate dans $ \Gamma_ \omega \backslash M = \Gamma_ \omega \backslash G/H\, .
$ \\

\noindent{\bf 5ème étape.} Dans toutes les opérations effectuées au cours des étapes
précédentes on sait que les données $ \omega ,L,J_o,J_oL $ sont invariantes par $
\Gamma_ \omega \, . $ Il reste à vérifier que la connexion linéaire $ D $ déterminée
par $ (L,J_oL ) $ est invariante par $ \Gamma_ \omega \, . $ \par \noindent Soit $
\Gamma \in \Gamma_ \omega \, . $ Posons $ D^\Gamma_ XY = \Gamma^{ -1}_\ast
D_{\Gamma_ \ast X}\Gamma_ \ast \, . $ On a alors

$$ \begin{array}{rl} (D^\Gamma_ X\omega )(Y,Z) & = X\omega( Y,Z) - \omega
(\Gamma^{ -1}_\ast D_{\Gamma_ \ast X}\Gamma_ \ast Y,Z ) - \omega (Y,\Gamma^{
-1}_\ast D_{\Gamma_ \ast X}\Gamma_ \ast Z ) \\  &  = X\omega( Y,Z) - \omega (\Gamma_
\ast Y,D_{\Gamma_ \ast X}\Gamma_ \ast Z ) - \omega (D_{\Gamma_ \ast X}\Gamma_ \ast
Y,\Gamma_ \ast Z ) \\  &  = 0\, .
\end{array}
$$

D'un autre c\^oté $ D^\Gamma L \subset  L $ et $ D^\Gamma J_oL \subset J_o\, L\, . $
\par \noindent Puisque $ D^\Gamma $ est aussi sans torsion on a $ D^\Gamma = D $ par
unicité.

\noindent{\bf Conclusion finale. }Les données $ (M,\omega ,L,J_oL,D ) $ sont
invariantes par l'action de $ \Gamma_ \omega $ (dans $ M)\, ; $ ces données donnent
lieu à une structure Kählerienne affinement plate dans $ \Gamma_ \omega \backslash
M\, ; $ d'o\`u le théorème 4.1.1.

\bigskip
(1) Pour établir la complète intégrabilité de $ {\mathcal N° }= J_oL $ on peut
éviter le recours au théorème de Bruhat--Whitney gr\^ace à l'analogue kählérien de
la réduction symplectique de Marsden--Weinstein, voir les sous-paragraphes
 ci-dessous et particulièrement le théorème de réduction Kählérien .

\subsection{Filtration symplectique de $(M,\omega)$}

Dans ce \S\ $ (M,\omega) $ est une variété symplectique connexe et $ G $-homogène
o\`u $ G $ est un groupe de Lie complètement résoluble.

\par \noindent On fixe une identification de $ (M,\omega) $ avec $ (G,H,\omega) $
comme dans les \S\ précédents. On suppose que $ \omega  \in  \Omega^ 2_l( G) $
possède un diagramme pondéré simple $ d(\omega ,F)\, . $

$$ \circ \rightarrow \circ \leftrightarrows \circ$$

Il est donc associé à $ (\omega ,F) $ la suite $ (G_k,H_k,\omega_ k ) $ ainsi que
les variétés symplectiques homogènes $ (M_k,\omega_ k )\, . $ On sait que $
(M_k,\omega_ k ) $ est déduite de $ (M_{k+1},\omega_{ k+1} ) $ par la méthode de
réduction de phase de MARSDEN--WEINSTEIN $\lbrack$MW]. On a le diagramme suivant

$$ \begin{matrix} {\mathcal M}_{k+1} &   \hookrightarrow &
  M_{k+1} \\  \downarrow &    &
\\  M_k &    &
\end{matrix}
$$

gr\^ace auquel on sait que $ M_k $ possède un feuilletage symplectique $ {\mathcal
D} $ de codimension 2 subordonné à l'hypersurface $ {\mathcal M}_{k+1}, $ i.e les
feuilles de $ {\mathcal D} $ qui rencontrent $ {\mathcal M}_{k+1} $ sont contenues
tout entières dans $ {\mathcal M}_{k+1}\, . $ De plus chaque feuille de $ {\mathcal
D} $ contenue dans $ {\mathcal M}_{k+1} $ est un revêtement symplectique de $ M_k\,
. $ On sait en outre que ces revêtements sont des revêtements des structures
bilagrangiennes affinement plates. Les variétés en jeu sont connexes et simplement
connexes, de sorte que l'espace des feuilles $ M_{k+1}/{\mathcal D} $ est une
surface symplectique qui hérite naturellement d'une structure bilagrangienne
affinement plate dès que $ M_{k+1} $ en possède une qui prolonge celles des feuilles
de $ {\mathcal D}\, . $
\par \noindent Les remarques ci-dessus permettent au moyen d'arguments de
récurrence de construire dans les $ (M_{k+1},\omega_{ k+1} ) $ des structures
Kähleriennes affinement plates. \par \noindent Ce procédé est utilisé dans [NB 4].

\section{Réductions Kählériennes}

\subsection{Digrammes normaux}

Soit $ (\omega  , F ) \in \Omega^2_l ( G ) x  {\mathcal F}( G )$ . Le diagramme
$d(\omega,F)$ est dit normal si en chaque sommet de coordonnées $(G_k,H_k)$ $H_k$
est un sous-groupe distingué de $G_k$ . Dans cette section on s'intéresse uniquement
à celles des 2-formes fermées qui sont de rang maximal. Pour le traitement dynamique
de ces formes on peut éviter le recours au théorème de BRUHAT--WHITNEY grâce à
l'analogue kählérien de théorème de réduction symplectique de MARSDEN--WEINSTEIN.

\par Soient $ G $ un groupe de Lie complètement résoluble de dimension paire $ 2m +
2 $ et $ \omega \in \Omega^ 2_l( G) $ une 2-forme de rang maximal. Si $ \omega $ est
cocompacte alors tout $ \Gamma_ \omega \in  {\mathcal C}_\omega( G) $ qui est
discret est un réseau cocompact dans $ G\, . $ Ci-dessus est un résultat
élémentaire.

\begin{lemma}
\label{lem:311} Soient $ G $  un groupe de Lie complètement résoluble et $ \omega $
une 2-forme fermée invariante à gauche de rang maximal. Alors $ \omega $  possède un
graphe orienté simple $ gr(\omega ,F) $  tel que pour chaque indice $ k\, , $  le
sommet $ H_k $  est distingué dans $ G_k\, . $
\end{lemma}

{\it Preuve}. On utilise un argument de récurrence élémentaire. La conclusion du
Lemme est trivialement vraie en dimension 2. En dimension $ 2m+2 $ on choisit dans $
G $ un sous-groupe connexe de dimension $ 1,G_1, $ qui est distingué dans $ G\, . $
Soit $ h_1 $ la sous-algèbre de Lie de l'algèbre de Lie de $ G $ qui est tangente à
$ G_1. $ On fixe un générateur $ \zeta $ de $ g_1 $ et on désigne par $ G_{2m+1} $
le sous-groupe de Lie connexe de $ G $ qui est tangent en l'élément neutre au noyau
de la 1-forme $ i(\zeta) \omega \, . $ Soit $ \underline{G}_{2m} $ le groupe de Lie
quotient de $ G_{2m+1} $ par $ G_1\, . $ Alors $ \underline{G}_{2m} $ hérite de $
\omega $ d'une 2-forme fermée $ \underline{\omega}_{ 2m} $ qui est de rang maximal.
L'argument de récurrence assure l'existence dans $ \underline{G}_{2m} $ d'un drapeau
en sous-groupes $ \underline{F} $ donnant lieu au graphe simple $ gr (
\underline{\omega}_{ 2m},F ) $ qui satisfait la conclusion du Lemme 2.1.1. On note $
\pi $ l'homomorphisme canonique de $ G_{2m+1} $ sur $ \underline{G}_{2m}\, ; $ alors
$ \pi^{ -1}( \underline{F}) \dot \subset  G $ a les propriétés requises . En effet
en vertu de l'hypothèse de récrrence  la série de composition $F := \underline{G}_l
\subset \underline{G}_2 \subset ... \subset \underline{G}_{2m}$ est telle que pour
chaque indice $k = 1 ,.. , 2m$ le sous-groupe $H_k$ est distingué dans $G_k$. Il en
résulte de l'image inverse par $\pi$ de $\underline{H}_k$ est distinguée dans celle
de $\underline{G}_k$. \\

\noindent{\it Remarque.} Conformément à la configuration (\ref{equation:8}) si le
groupe $ G $ est nilpotent, alors chaque 2-forme cocompacte $ \omega $ de rang
maximal définie dans $ G $ donne lieu à une 2-forme cocompacte $
\underline{\omega}_{ 2m} $ dans le groupe $ \underline{G}_{2m}\, . $ Il suffit pour
cela que $ \Gamma_ \omega \cap  G_1 $ soit uniforme dans $ G_1 $ convenablement
choisi dans $G$ . En effet puisque $G$ est nilpotent il existe un sous-groupe à un
paramètre central dans chaque qui
intersecte le réseau donné suivant un réseau dans lui même. \\

\noindent Dans la suite on conserve un graphe $ gr(\omega ,F) $ qui satisfait le
Lemme \ref{lem:311} ci-dessus. On a ainsi le diagramme suivant

\begin{equation}
\label{equation:8}
\begin{matrix}  G_1 &    &   G_1 \\
\downarrow  &    &  \downarrow \\
G_{2m+1} &   \stackrel{i}{\rightarrow} & G \\
\downarrow \pi &    & \downarrow \\
\underline{G}_{2m} & \rightarrow & \underline{G}_{2m+1}
\end{matrix}
\end{equation}

et la relation $ i^\ast \omega  = \pi^ \ast \underline{\omega}_{ 2m}\, . $ On se
propose de mettre en évidence l'analogue k\"ahlérien de la réduction symplectique de
MARSDEN--WEINSTEIN.

 A partir de maintenant la situation est la suivante .
On fixe un groupe de Lie complètement résoluble $G$ qui est connexe et simplement
connexe . On suppose que    $\omega \in \Omega^2_l( G )$  est une forme symplectique
cocompacte . On fixe une fois pour toutes un réseau  $\Gamma_{\omega}$   dans
$C_{\omega}(G)$ , un diagramme normal $d(\omega ,F)$ ainsi qu'un feuilletage
lagrangien $L$ qui est invariant par les translations à gauche dans $G$ . On fixe
une structure bilagrangienne affinement plate $( G , L  , N = J(L) )$ . Soit  $h =
g_{J_0} + \sqrt{-1}\omega$ une métrique de Kähler invariante par les translations à
gauche dans $G$ définies par les éléments du réseau $\Gamma_{\omega}$. Conformément
aux notations de la Section 1 soit $J_o$ le tenseur de la structure presque complexe
de $h$ . On se propose de démontrer la complète intégrabilité de $J_o(L )$ sans
recourir au théorème de Bruhat-Whitney. \\

\noindent Notons $ {\mathcal T}_1 $ et $ {\mathcal T}_2 $ les sous-fibrés vectoriels
de $ TG $ engendrés par $ \tilde \zeta $ et par $ J_o\tilde  \zeta $ respectivement.
On note $ {\mathcal T}_o $ le sous-fibré vectoriel orthogonal à $ {\mathcal
T}_1\oplus  {\mathcal T}_2 $ relativement à la métrique riemannienne $g_{J_0}.$ On a
la somme orthogonale

$$ TG = {\mathcal T}_o\oplus  {\mathcal T}_1\oplus  {\mathcal T}_2\, . $$

En tout $ \Gamma  \in  G $ les sous-espaces vectoriels $ {\mathcal T}_o(\Gamma) $ et
$ {\mathcal T}_1(\Gamma) $ sont tangents à la sous-variété $ \Gamma \cdot G_{2m+1}\,
. $ En particulier, en restriction à $ G_{2m+1}, $ $ {\mathcal T}_1 $ n'est pas
autre chose que le fibré vertical de la fibration principale

$$ G_1\rightarrow  G_{2m+1}\rightarrow  \underline{G}_{2m} $$

et $ {\mathcal T}_o $ en est le fibré orthogonal pour la métrique riemannienne
induite par $ g_J\, . $ Naturellement $ {\mathcal T}_o $ est un sous-fibré du fibré
vectoriel complexe $ TG\, . $ Si $ s $ est une section de $ {\mathcal T}_o $ alors
on a $ \omega (\tilde \zeta ,s ) = \tilde  \zeta (J_o\tilde \zeta ,s ) = 0\, . $

\subsection{Quelques résultats auxiliaires}

 L'objet principal de cette sous-section est de démontrer sans recourir au
théorème de Bruhat-Whitney que les systèmes différentiels définis respectivement par
$T_o$ et $T_l + T_2$ sont complètement intégrables . Comme indiqué plus haut on
gardera en mémoire la configuration (\ref{equation:8}) reproduit ci-dessous.

$$
\begin{array}{ccc}
                   G_1    &  &   G_1 \\
                   \downarrow & & \downarrow \\
                   G_{2m+1}  & \hookrightarrow &  G \\
                    \downarrow & & \downarrow \\
                   G_{2m} & \hookrightarrow  &    G_{2m+1}
\end{array}
$$

Avec $\dim(G) = 2m + 2$ . L'action de $G_1$ par les translations à gauche dans G est
hamiltonienne .

La variété symplectique $( G_{2m} ,\omega )$ est une réduite de $( G ,\omega  )$
pour cette action . Nous rappelons que $G$ est muni des structures Kählériennes $( G
,\omega , J )$ et $( G ,\omega   , J_o )$ . Les variéts symplectiques $(G ,\omega )$
et $(G_{2m} ,\omega )$ sont munies de des structures bilangrangiennes affinement
plates $( G ,\omega   , L , J(L) )$ et $( G_{2m} ,\omega   , L , J(L) )$
respectivement.

On se place désormais dans un ouvert connexe U  domaine des coordonnées de
Darboux-Hess  $( q_1 ,..,q_{m+1} , p_1 , .. ,p_{m+1} )$.
 En d'autres termes les fonctions  $q_i$ et $p_i$ satisfont les
conditions suivantes  (10) \\
       (i) $ \omega   =  \sum_i dq_i  \wedge dp_i$ , \\
       (ii)  les $p_i$ sont des intégrales premières du feuilletage langrangien $L$, \\
       (iii)  les fonctions $q_i$ sont des intégrales premières de $J(L)$ , \\
       (iv)  $(q_i ,p_i )$ est un système des coordonnées locales affines de la
 structure affine définie par $(L ,J(L) )$ ; de plus si $( q'_i , p'_i )$ est
 un autre système de coordonnées locales de Darboux-Hess défini dans le même
domaine $U$ alors le changement de coordonnées $q' = q'(q,p)$ , $p'=p'(q,p)$ est une
transformation affine rigide de la structure bilangrangienne affinement plate $( L ,
J(L) ) $.

Une fois pour toutes on fixe un système des coordonnées  de Darboux-Hess adapté à $(
L , J(L), (q_i,p_i)$ ; on désigne par $S_i$ le champ de vecteurs hamiltonien de la
fonction $ p_i$ . La propriété (i) ci-dessus assure que $J(S_i)$ est le champ de
vecteurs hamiltonien de la fonction $q_i$ . Il en resulte que le système des champs
de vcecteurs $(S_i , J(S_i ) )$ est orthonormée pour la métrique riemanienne $g(u,v)
= ( Ju ,v )$ . Nous allons maintenant mettre en évidence des propriétés
intéressantes de champs de vecteurs $S_i$.

\begin{lemma}
\label{lem:321} Les champs de vecteurs $S_i$ et le tenseur $J_o$ sont liés comme il
suit (a) $[ S_i , J_o S_j ]$ est une section de $J_o(L)$, (b) $[ J_o S_i , J_o S_j
]$ est une section de $L$ .
\end{lemma}

Démonstration . (a) On désigne par $R_u$ le fibré vectoriel trivial $U \times
{\mathbb R}^{2m+2}$ et par $L_u$ la restriction à l'ouvert $U$ de $L$ . Soit
$\theta_o$ l'homorphisme de fibré vectoriel de $L_u$ dans $R_u$ défini par :

$$\theta (x,v ) = (g_o (S_1 , v ) ,.., g_o (S_{m+1} , v ) , g_o (JS_1 , v ),.., g_o (
JS_{m+1} , v )).$$

On rappelle que $g_o (S_i ,v ) = \omega (J_o S_i , v )$ . L'application $\theta$ est
un plongement de $L$ dans $R_u$ , par conséquent l'image de Lu est un sous-fibré
vectoriel trivial de $R_u$ . On choisit une base des sections locales constantes de
$\theta(L_u )$, soit $c_1 ,.. , c_{m+1}$ ; en d'autres termes chaque $c_j$ est une
application constante de l'ouvert $U$ dans ${\mathbb R}^{2m+2}$. On désigne par
$\alpha_j$ l'image inverse par $\theta$ de la section $c_j$. On obtient ainsi une
nouvelle base des sections de $L_u$ satisfaisant les conditions

\begin{equation}
\label{equation:9}
 \omega( J_o S_i , \alpha_j ) = constante \ \ \ and \ \ \
 \omega( J_o JS_i , \alpha_j ) = constante.
\end{equation}

On utilise maintenant la métrique Riemanienne $g(u,v) = \omega(Ju,v)$ et sa base *
orthonormée $(S_i , JS_i )$ pour décomposer les $(\alpha)_i$ et les $J_o (\alpha_j
)$:

$$J_o (\alpha_i )=\sum_k \omega(JS_k ,J_0 \alpha_i )S_k + \sum_k \omega(J_o
\alpha_i ,S_k )JS_k. $$

Compte tenu des conditions (\ref{equation:9}) ci-dessus on a $[
S_k , J_o \alpha_j ] = O$ pour $j , k = 1 ,.. , m+1 $. Maintenant
considérons un triplet $(\alpha_i , S_j , S_k )$ ; on déduit de la
nullité de $[ S_j ,J_o \alpha_k ]$ l'identité suivante
$\omega([S_j,J_o \alpha_i] , J_o S_k ) = \omega ([ S_j ,J_o S_k ]
,J_o \alpha_i ) = o $. Puisque les $J_o \alpha_i$ engendrent le
sous-fibré lagrangien $J_o (L)$ on en déduit que les $[ S_j ,J_o
S_k ]$ sont des sections de $J_o (L)$ , (a) est démontré.

(b) La seconde assertion découle directement de (a) ; en fait puisque les champs de
vecteurs $S_j$ commutent deux à deux on a l'identité suivante $[J_o S_i , J_o S_j] =
J_o ( [J_o S_i , S_j ] + [ S_i , J_o S_j ] )$ d'où il découle que les $[ J_o S_i ,
J_o S_j ]$ sont des sections de $L$. Cela termine la démonstration du Lemme
\ref{lem:321}.

Nous sommes maintenant en mesure de démontrer que $T_o$ est complètement intégrable.
Nous supposons que le champ de vecteur $S_1$ est une base de section de $T_1$ . On
utilise de nouveau la métrique Riemannienne $g_o$ pour décomposer $L$ comme il suit

\begin{equation}
\label{equation:10}
L= T_1 + orth_o ( T_1 ),
\end{equation}

où $orth_o (T_l )$ est le sous-fibré orthogonal de $S_1$ relativement à la métrique
$g_o$ . Considérons l'homomorphisme de fibré vetoriel $\phi$ de $L_u$ dans le fibré
trivial $U \times {\mathbb R}^{m+1}$ défini comme il suit : $\phi (tS_1 ,v ) = ( g
(tS_1 , v ) + tg_o (S_1 ,S_1 ),g(S_2 ,v ) +..+ g(S_m+1 ,v ))$ où $(tS_1 ,v )$ est
une section de $T_1 \times orth_o (T_1 )$ . L'application $\phi$ est un isomorphisme
de fibré vectoriel . L'image de $orth_o (T_1 )$ est un sous-fibré vectoriel de rang
$m$ de $U \times {\mathbb R}^{m+1}$ . Quitte à réduire U à un domaine de
trivialisation de $\phi ( orth_o (T_1 )$ on choisit une base de section de $\phi (
orth_o (T_1 )$ formée des sections constantes ,$w_1 ,.. , w_m$  ; chaque $w_j$ est
une application constante de $U$ dans ${\mathbb R}^{m+1}$ . Soit $\beta_j$ l'image
inverse par $\phi$ de $w_j$ . Compte tenu de (\ref{equation:10}) les $\beta_j$  sont
des sections de $( O ) \times orth_o (T_1 )$ , ces sections satisfont les conditions
suivantes $g (S_k ,\beta_j ) = constante$ , pour $j := 1 ,.., m$ et $k:= 1 ,.. ,
m+1$ . Dans la base $g$-orthonormée $( S_i , JS_i ) , i:= 1 , ..,m+1$ on décompose
les $\beta_j = \sum_k \omega (JS_k ,\beta_j )S_k $. Puisque les $\omega ( JS_k ,
\beta_j )$ sont sont constantes on a l'identité

\begin{equation}
\label{equation:11} [ S_k ,\beta_j ] = 0 , j:= 1,..,m ,k:= 1 ,.., m+1.
\end{equation}

 Retournons maintenant à la configuration (\ref{equation:8}) . Le système $( S_1 ,..,S_m ,JS_1,
..,JS_m )$ engendre le sous-fibré $g$-orthogonal de $T_1$ dans $TG_{2m+1}$ .

\begin{thm}
\label{thm:322} Les systèmes différentiels définis respectivement par  $T_1 + T_2$
et par $T_o$ sont complètement intégrables .
\end{thm}

Démonstration . Nous considérons la variété holomorphes $(G , J_0 )$ munie de la
métrique Kählérienne $h = g_o +  \sqrt{-1} \omega$ . Les sous-fibrés vectoriels $T_1
+ T_2$ et $T_0$ sont analytiques complexes . En vertu d'un Théorème de David L.
Johnson ,[JD] ,la complète intégrabilité de $T_0$ équivaut à la complète
intégrabilité de $T_1 + T_2$ . On va montrer que $T_1 + T_2$ est complètement
intégrable . On considère la base $g$-orthonormée $( S_i , JS_i  )$ et la base
$(\beta_j)$ de $orth_0 (T_1 )$ construite ci-dessus . On déduit de l'identité
(\ref{equation:11}) ci-dessus et du Lemme \ref{lem:321} les identités suivantes
$\omega ([S_1 , J_0 S_1 ] , \beta_j ) = 0$ et $ \omega ( [S_1 , J_0 S_1 ] , J_0
(\beta_j )) = 0$ . Cela montre que le crochet $[ S_1 ,J_0 S_1 ]$ est une section du
sous-fibré engendré par le système $(\beta_j , J_0 (\beta_j )$ qui n'est pas autre
chose que $T_0$ . Puisque $T_0$ et $T_1 + T_2$ sont $g_0$ -orthogonaux le crochet
$[S_1 , J_0 S_1]$ est une section de $T_1 + T_2$ . Le théorème de David Johnson
évoqué ci-dessus assure que $T_0$ et $T_1 + T_2$ sont complètement intégrables , ce
qui termine la démonstration du Théorème \ref{thm:322}. \\

On va se placer maintenant dans la variété analytique complexe $( G , J_0 )$ . Elle
est équipée de deux feuilletages holomorphes transverses $T_1 + T_2$  et $T_0$ . En
outre le feuilletage $T_0$ est subordonné au sous-groupe de Lie $G_2m+1$ dans le
sens que toute sous-variété intégrale maximale de $T_0$ contenant un point de
$G_{2m+1}$ est entièrement contenue dans $G_{2m+1}$ . Les feuilles de $T_0$
contenues dans $G_{2m+1}$ sont des revêtement symplectiques de la variété réduite $(
G_2m , \omega )$ de la configuration (\ref{equation:8}) . En fait le sous-fibré
$T_0$ est une connexion principale de la $G_1$ -réduction   $G_{2m+1}
\longrightarrow G_{2m}$ qui est une $G_1$ - fibration principale . Puisque le groupe
de Lie $G$ est simplement connexe et complètement resoluble il en est de même de
$G_{2m}$  qui hérite ainsi de chaque de $T_0$ d'une métique Kählérienne . Etant
donné que $T_0$ est une $G_1$ -connexion principale plate la structure Kählérienne
de $G_{2m}$ héritée d'une nappe d'holonomie de $T_o$ est indépendante de la nappe
choisie . On peut donc énoncer ce qui suit .

\begin{thm}
\label{thm:323} La configuration (\ref{equation:8}) est une réduction Kählérienne de
la variété Kählérienne $(G,\omega,J_0)$ .
\end{thm}

{\it Remarque}. On sait déjà par construction de la métrique Kählérienne $h_ = g_ +
\sqrt{-1} \omega$ que cette dernière possède une réduction dans $G_{2m}$ . Ainsi les
deux métriques $h$ et $h_0$ possèdent des réductions dans $G_{2m}$ . Nous continuons
de noter ces deux structures $( G_{2m} , \omega , J )$ et $( G_{2m} ,\omega , J_0 )$
respectivement . Le sous-fibré lagrangien de $TG_{2m}$ hérité du sous-fibré $L$ de
$TG$ est aussi  abusivement noté $L$ .

Nous allons maintenant démontrer sans nous servir du théorème de
Bruhat-Witney le resultat suivant.

\begin{thm}
\label{thm:324} Le sous-fibré lagrangien $J_0 (L)$ est complètement intégrable .
\end{thm}

{\it Démonstration}. Compte tenu du théorème \ref{thm:323} on va procéder par
récurrence sur la dimension du groupe de Lie G . En dimension 2 il n'y a rien à
démontrer. Suppososns maintenant le théorème \ref{thm:324} vrai en dimensions
inférieures à 2m+2 .Il est alors appliquable dans $( G_{2m} , \omega , J_0 (L) )$ .
On supposera alors sans perte de généralité que les coordonnées de Darboux-Hess
$(q_i , p_i ) ,i = 1 , .. , m+1$ sont choisies de sorte que les fonctions  $q_2
,..,q_{m+1} ,p_2 ,.., p_{m+1}$ se projettent en un système des coordonnées de
Darboux-Hess pour $( G_{2m} , \omega , L ,J(L )$ . Puisque Le sous-fibré $J_0 (L)$
de $TG_2m$ est complètement intégrable il résulte
 du Lemme \ref{lem:321} que les champs de vecteurs $J_0 S_i$ pour $i = 2
,..,m+1$ commutent deux à deux . Le même Lemme \ref{lem:321} assure que pour tout $i
= 2 ,.., m+1$ les crochets $[S_1 , J_0 S_i ]$ et $[S_i , J_0 S_1 ]$ sont des
sections de $J_0 (L)$ et que $[ J_0 S_1 , J_0 S_i]$ est une section de $L$ . La
métrique Riemannienne $g_0$ étant  quasi fibrée les feuilletages $T_0$ et $T_1 +
T_2$ sont totalement géodésiques par rapport à la connexion de Levi-Civita de $g_0$
,( [JW ] ). En combinant le lemme 2.3.1 et les propriétés géodésiques rappelées
ci-dessus on déduit que les crochets $[ J_0 S_1 ,J_0 S_i]$ sont nuls . Cela termine
la démonstration du théorème \ref{thm:324}.

Maintenant nous sommes en présence d'une structure bilagtrangienne dans $( G ,
\omega )$ définie par la paire $( L , J_o (L))$ des feuilletages lagrangiens dans $(
$ G $ , \omega)$. A la section 1 on a démontré par un calcul direct la nullité du
tenseur de courbure de la connexion linéaire définie par la paire $( L,J_0(L))$ . La
procédure de réduction Kählérienne de la section 2 contient en fait des ingrédients
qui permettent de démontrer rapidement que la structure bilangangienne $( G , \omega
, L , J_0 (L) )$ est affinement plate ; nous nous proposons d'en exhiber des
coordonnées de Darboux-Hess .En fait considérons les coordonnées de Darboux-Hess de
$( G , \omega , L , J(L) )$ déjà fixées  $( q_i ,p_i ) , i = 1 ,.., m+1$  ainsi que
les champs de vecteurs hamiltonienns $S_i$ des fonctions $p_i$ . Le Lemme
\ref{lem:321} joint au Théorème \ref{thm:324} a permis de voir que les champs $S_i$
commutent deux à deux et que les crochets $[ S_i ,J_0 S_j ]$ sont des sections de
$J_0 (L )$ . Considérons de nouveau l'isomporphisme de fibré vectoriel de $L_u$ sur
le fibré trivial $U \times {\mathbb R}^{m+1}$ défini par $ \theta ( v ) = ( g_0 (S_1
, v ),.. , g_0 (S_{m+1} , v )$ où $v$ est une section de $L$ . Soit $c_j$ la section
constante de $U \times {\mathbb R}^{m+1}$ qui envoie tout point de $U$ sur le jième
vecteur de la base canonique de ${\mathbb R}^{m+1}$ . On obtient ainsi une base des
sections de $U \times {\mathbb R}^{m+1}$ . Soit $\alpha_j$ l'image inverse par
$\theta$ de la section $c_j$ ; les sections $\alpha_j$ de $L_u$ vérifient les
conditions suivantes $ \omega (J_0 S_j , \alpha_k ) = \delta_{jk} , j ,k = 1 ,.. ,
m+1 $.

\begin{prop}
\label{prop:325} Les champs de vecteurs $J_0 (\alpha_j )$ sont localement
hamiltoniens.
\end{prop}

{\it Démonstration}. En vertu du lemme \ref{lem:321}  on vérifie par un calcul
direct facile que la différentielle extérieure de $\omega (J_0\alpha_i , . )$ est
nulle sur les couples $(S_i ,S_k ) , (S_i , J_0 S_k )$ et $(J_0S_i , J_0 k )$ . La
proposition \ref{prop:325} est démontrée .

On déduit de la proposition \ref{prop:325} le résultat attendu, à savoir le théorème
qui suit.

\begin{thm}
\label{thm:326} Les fonctions numériques locales définies par le système des champs
de vecteurs hamiltoniens $( S_j , J_0 \alpha_j )$ sont des systèmes de coordonnées
de Darboux-Hess de la structure bilagrangienne $( G , \omega , L ,J_0 (L) )$.
\end{thm}

{\it Démonstration}. En fait on a les identités suivantes
$$\omega ( S_i , [S_j ,J_0 \alpha_k ] ) = S_j \omega ( S_i ,J_0 \alpha_k ) = 0.$$
On en déduit que les crochets
$[ S_j , J_0 \alpha_k ]$ sont des sections de $L $. En vertu du Lemme \ref{lem:321}
$[ S_j , J_0 S_k ]$ est une section de $J_0 (L)$ , on a donc les identités $\omega
([ S_j , J_0 \alpha_i ] , J_0 S_k ) = \omega ([ S_j , J_0 S_k ] , J_0 \alpha_i ) = 0
$. Cela montre que les crochets $[ S_j , J_0 \alpha_k ]$ sont tous nuls , $i ,j = 1
,.., m+1$ . La proposition \ref{prop:325} entraine que les champs de vecteur  $J_0
\alpha_i$ commutent deux à deux . Par conséquent le système $( S_i ,.., S_{m+l} ,J_0
\alpha_1 ,.., J_o \alpha_{m+l} )$ est une base symplectique de $( G ,\omega ))$. Par
construction des $S_i$ et des $\alpha_i$ les fonctions hamiltoniennes locales
définies par le système $( S_i , J_0 \alpha_i )$ forment un système des coordonnées
de Darboux-Hess de $( G ,\omega , L , J_0 (L))$ . Cela achève la démonstration du
théorème \ref{thm:326}.

On déduit du théorème \ref{thm:326} le corollaire important qui suit.

\begin{thm}
\label{thm:327} La structure bilagrangienne affinement plate $( G , \omega , L , J_0
(L) )$ determiné par le théorème \ref{thm:326} ci-dessus est invariante par les
translations à gauche dans $ G $ définies par les éléments du réseau $
\Gamma_{\omega} $ .
\end{thm}

Démonstration . Les feuilletages lagrangiens $L$ et $J_0 (L)$ sont invariant par
l'action à gauche dans $ G $ du réseau $\Gamma_{\omega} $. Il en est de même de
l'unique connexion sympectique sans torsion qui est associée à la paire $( L , J_0
(L))$.

\section{Groupes Kärleriens}

 Cette section 4 est consacré à des
applications des résultats obtenus dans les précédentes, en particulier ceux de la
section 3.

\subsection{Réseaux Kähleriens des groupes de Lie résolubles}

L'expression ``groupe kählerien'' est utilisée dans le sens de [ABCKT]. Il s'agit
des groupes $ \Gamma $ qui sont isomorphes aux groupes fondamentaux des variétés
kähleriennes compactes $ (M,\omega ,J)\, . $

\par \noindent On ne sait pas caractériser complètement les groupes
kähleriens. La question traitée ici est d'une nature autre. Elle est liée en fait à
la dynamique différentiable dans la variété k\"ahlérienne $ M\, . $
\par
\bigskip
\noindent{\bf Question 1.}  Quelles sont les variétés k\"ahlériennes compactes $
(M,\omega ,J) $  dont la variété analytique réelle sous-jacente $ M $  admet une
action localement libre transitive d'un groupe de Lie complètement résoluble.

 \noindent Nous allons répondre à la
question moyennant une hypothèse supplémentaire qualifiée d' Imperfection
géométrique. \par \noindent Soit $ G $ un groupe de Lie solution de la question 1
ci-dessus ; soit $ g $ l'algèbre de Lie de $ G\, . $ En vertu de l'identification de
l'algèbre de de RHAM $ H^\ast_{ DR}(M,\bkR) $ avec celle de KOSZUL $ H^\ast( g,\bkR)
, $ ($\lbrack$HA],$\lbrack$RA]) la classe de K\"ahler $ [\omega] \in
H^2_{DR}(M,\bkR) $ de $ (M,\omega ,J) $ contient des formes homogènes $ \omega_ o\in
\Omega^ 2_l( G)\, . $ L'imperfection géométrique consiste à supposer que la forme de
K\"ahler $ \omega $ est $ G $-homogène. Rappelons de nouveau la signification de
formes symplectque homogène définie dans un espace homogène, ([McD1]). Soit un
groupe de Lie complètement résoluble contenant un réseau cocompact  $\Gamma$. On
fait opérer le réseau $\Gamma$ dans $ G $ par translation à gauche . L'espace des
orbites de cette action de $\Gamma$ est un $ G $-espace homogène à droite . Si
$\omega$ est une forme symplectique dans  $\Gamma \backslash G $ alors l'imge
inverse par la projection canonique de $ G $ sur $\Gamma\backslash G $ est une forme
symplectique invariant par les translations à gauche définies par les éléments du
réseau $\Gamma$, on dira que la forme symplectique $\omega$ de $\Gamma \backslash G
$ est $G$-homogène si son image inverse par $\pi$ est invariante par toutes les
translations à gauche dans le groupe $ G $ .

\noindent Dans la mesure où le type d'homotopie rationnelle de $ M $ ne dépend que
de l'algèbre de de RHAM, l'imperfection géométrique est topologiquement neutre, pour
davantage de détails le lecteur consultera [BG 1], [BG 2], [McD1].

\noindent Du point de vue des groupes k\"ahleriens la question 1 a des variantes. En
voici une légèrement (en apparence !) plus faible. \par
\bigskip
\noindent{\bf Question 2.}  Quelles sont les variétés K\"ahleriennes compactes $
(M,\omega ,J) $  dont les groupes fondamentaux sont isomorphes à des réseaux
uniformes dans des groupes de Lie complètement résolubles . Les exemples suivants
sont une illustration des questions 1 et 2 ci-dessus . Considérons l'espace
vectoriel complexe ${\mathbb C}^2$ muni du produit hermtien canonique . On considère
maintenant l'ensemble ${\mathbb Z}^4 $ muni de la loi de composition suivante : $(
p_1 , p_2 ,p_3, p_4 ). (q_1 , q_2 , q_3 , q_4 ) = (p_1 ,p_2 ,p_3 , p_4 ) +
 exp( \sqrt{-1} p_4 \pi/2)(q_1 ,q_2 ), q_3 , q_4 ) .$
On obtient ainsi un groupe discret qui agit proprement et holomorphiquement dans
${\mathbb C}^2$ comme il suit : $(p_1 ,p_2 ,p_3 ,p_4 ). (z_1 ,z_2 ) = (p_1 ,p_2 ,p_3
,p_4 ) + ( z_1 exp( \sqrt{-1}p_4 \pi/2 ) , z_2 )$ . Le produit hermitien de
${\mathbb C}^2$ est invariant par l'action de ${\mathbb Z}^4$ . Ainsi l'espace des
orbites obtenu est une variété Kählérienne compacte qui n'admet pas d'action
localement simplement transitive de groupe de Lie complètement résoluble . Pour s'en
convaincre il suffit de calculer le premier nomnbre de Betti de cette variété qui
est inférieur à 4 . D'après [ HA ] (voire aussi [ RA ] ) si une telle action
existait alors le premier nombre de Betti de l'espace des orbites vaudrait 4 . \par
\bigskip
\noindent Voici un résultat en rapport avec les questions 1 et 2.

\begin{thm}
\label{thm:401}
  Soit $ \Gamma $  un réseau cocompact dans le groupe de Lie complètement
résoluble $ G $  et soit $ \omega \in \Omega^ 2_l( G) $ une forme symplectique. Si
la variété symplectique $ (\Gamma \backslash G,\omega) $ possède une métrique
k\"ahlérienne $ (M,\omega ,J) $  alors $ \Gamma $  est virtuellement commutatif.
\end{thm}

\noindent En vertu de l'arsenal mis au point dans les \S\ précédents on déduit le
théorème \ref{thm:401} du suivant.

\begin{thm}
\label{thm:402}
 Soient $ \Gamma $  un réseau cocompact dans le groupe de Lie complètement
résoluble $ G $  et $ \omega  \in \Omega^ 2_l( G) $
 une forme symplectique. Si la variété symplectique $ (\Gamma \backslash
G,\omega) $  possède une structure k\"ahlerienne $ (\Gamma \backslash G,\omega ,J_0
) $  alors la variété réelle sous-jacente $ \Gamma \backslash G $  est difféomorphe
à un quotient fini du tore plat
\end{thm}

 \noindent Avant de tirer d'autres conséquences des théorèmes  \ref{thm:401}
 et \ref{thm:402},
nous allons en donner les démonstrations. \\

{\it Démonstration du théorème \ref{thm:402}}. D'après nos hypothèses la forme de
K\"ahler $ \omega $ de $ (\Gamma \backslash G,\omega ,J_o ) $ est $ G $-homogène. En
vertu du théorème 2 de [BG2 ] la 2-forme $ \omega $ possède un diagramme pondéré
simple

$$\circ \rightarrow \circ \leftrightarrows \circ $$

La variété symplectique homogène $ (G,\omega) $ possède un feuilletage lagrangien $
L\, , $ invariant par les translations à gauche dans $ G\, . $  Nous pouvons
appliquer les ingrédients dévélopés à la Section 3. Pour cela nous adoptons les
notations des théorèmes \ref{thm:326} et \ref{thm:327}. En vertu du théorème
\ref{thm:327}, $ (G,\omega) $ possède une structure Kählerienne affinement plate
$(G,\omega ,J ,D)\, ; $ les feuilles de $ L $ étant des sous-variétés totalement
réelles de $ (G,J ) $ et totalement géodésiques de $ (G,D )\, . $ Naturellement ni $
J\, , $ ni $ D $ ne sont invariants par $ \Gamma \, . $ L'ensemble $ {\mathcal
K}^a_\omega( G) $ est donc non vide. Par ailleurs $ (G,\omega) $ hérite de $ (\Gamma
\backslash G,\omega ,J_o ) $ d'une structure k\"ahlerienne notée aussi $ (G,\omega
,J_o ) $ qui est $ \Gamma $-invariante. A ce point on sait d'après le théorème
\ref{thm:231} que l'intersection $ K_\omega( G) \cap K^a_\omega( G) $ est non vide.
En utilisant auxilairement les structures Kählériennes  $(G,\omega,J )$ et
$(G,\omega ,J_0 )$ d'une part et la structure bilagrangienne affinement plate
$(G,\omega ,L, J(L) )$ d'autre part nous avons démontré que la paire $(L,J_0(L))$
définit dans $( G , \omega )$ une structure Kählérienne affinement plate voir
Théorème \ref{thm:327}. En d'autres termes, $ (G,\omega) $ possède une structure
k\"ahlerienne affinement plate $ (G,\omega ,\hat J ) $ qui est $ \Gamma $-invariante
(Théorèmes \ref{thm:326} et \ref{thm:327}). Cette structure $(G , \omega , L , J_0
(L))$ donne naissance dans la variété compacte $ \Gamma \backslash G $ à une
métrique Kählérienne affinment plate [NB4].
\par

\smallskip
\noindent En vertu des théorèmes classiques de clasification des variétés des
variétés Riemanniennes plates , [WJ], $ \Gamma \backslash G $ est un quotient fini
du tore plat. Cela termine la démonstration du théorème \ref{thm:402}

\noindent En vertu du théorème \ref{thm:402} le groupe fondamental du tore plat est
isomorphe à un sous-groupe distingué commutatif $ \Gamma_ o\subset  \Gamma $ qui est
d'indice fini dans $ \Gamma $. Cela n'est pas autre chose que le contenu du théorème
\ref{thm:401}.

\noindent On observe que $ \Gamma_ o\backslash \Gamma $ est en fait isomorphe au
groupe d'holonomie linéaire de la variété affinement plate $ (\Gamma \backslash
G,D)\, . $ On a donc du théorème \ref{thm:402} d'autres corollaires remarquables
dont le suivant
\par
\bigskip

\begin{thm}
\label{thm:403} Soient $ \Gamma $  un réseau cocompact dans le groupe de Lie
complètement résoluble $ G $ et $ \omega  \in \Omega^ 2_l( G) $
 une forme symplectique. Si la variété symplectique $ (\Gamma \backslash
G,\omega) $  possède une métrique k\"ahlerienne, alors $ \Gamma $  est un groupe
cristallographique
\end{thm}

\noindent Ce théorème est une autre expression du théorème \ref{thm:402} , à savoir
que $ {\mathcal K}_\omega( G) \cap {\mathcal K}^a_\omega( G) $ est non vide.
\par \noindent Naturellement la connexion $ D $ est géodésiquement complète.
On peut déduire ce fait de la géométrie affine de $ (\Gamma \backslash G,D)\, , $
(voir $\lbrack$FGH]). De la relation $ D\omega  = 0 $ résulte que l'holonomie
linéaire de $ (\Gamma \backslash G,D) $ habite le groupe unimodulaire $ SL (\omega^
m )\, . $ Des classiques théorèmes de Bieberbach joints au théorème \ref{thm:403} et
au théorème de rigidité des réseaux [ON] entraînent les observations ci-dessous.

\par \noindent Soit $ G $ un groupe de Lie complètement résoluble, $ {\rm
Res}(G) $ désigne l'ensemble des réseaux uniformes $ \Gamma $ dans $ G $ et $ {\rm
Rig}(G) $ l'ensemble de leurs classes d'isomorphisme. Si deux réseaux $ \Gamma $ et
$ \Gamma ^\prime $ sont isomorphes alors $ \Gamma \backslash G $ possède une
métrique k\"ahlerienne si et seulement si $ \Gamma ^\prime \backslash G $ en possède
une. Cela résulte du théorème de rigidité des $ \Gamma $ dans $ G\, . $

Pour les propriétés classiques des groupes crystallographiques et les groupes de
Bieberbach on pourra consulter  [ WJ ] . Nous considérons les structures
k\"ahleriennes à formes de Kähler $ G $-homogènes. Puisque ces groupes sont
critallographiques, un des théorèmes fondamentaux de Bieberbach [W] assure ce qui
suit.

\begin{prop}
\label{prop:404}
 Le sous-ensemble de $ {\rm Rig}(G) $ constitué des classes des réseaux
k\"ahleriens dans $ G $ est fini.
\end{prop}

\subsection{Conjectures des tores plats} 

 Dans cette sous-section nous allons mettre en évidence quelques applications des études
menées dans les sections précédentes . Nous montrerons comment le théorème de C.
Benson et C. Gordon se déduit rapidement de notre théorème \ref{thm:323}. Auparavant
on commencera par s'intéresser à la conjecture des tores plats énoncée dans
l'introduction . Cette conjecture a motivé partiellement [ BN6 ] et le travail
présent . Pourvu des résultats mis au point dans les sections 1 , 2 et 3 nous sommes
maintenant en mésure de démontrer rapidement le théorème suivant

\begin{thm}
\label{thm:411} Soient $ \omega  \in  \Omega^ 2_l( G) $  une forme symplectique et $
\Gamma \subset G $ un réseau cocompact. Si la variété symplectique $ (\Gamma
\backslash G,\omega) $ possède une métrique k\"ahlerienne, $ (\Gamma \backslash
G,\omega ,J_0 )\, , $ alors la variété réelle $ \Gamma \backslash G $  est
difféomorphe au tore plat
\end{thm}

{\it Démonstration}. En fait la démonstration que nous allons donnée est une
combinaison de nos théorèmes \ref{thm:401} et \ref{thm:402} avec un  théorème de
HATTORI--RAGHUNATAN (et un théorème de Bochner (voir [AB]) si on préfère utiliser
l'argument de la dimension des ``champs de Killing''). En effet conformément aux
notations utilisées dans les sections 1 , 2 et 3 nous désignons encore par $( G ,
\omega , J_0 )$ la structure Kählérienne de dans $ G $ qui est $\Gamma$ -invariante
et qui se projette en la struture $( \Gamma \backslash G , \omega , J_0 )$ . Puisque
la forme de Kähler $\omega$ est $ G $ -homogène la machinérie mise au point dans les
sections 2 et 3 peuvent être utlisées pour construire une structure Kählérienne
affinement plate que nous désignons par $( G , \omega , \hat{J} )$ et qui est
définie par la structure bilagrangienne affinement plate $( G , \omega , L , J_0 (L)
)$  ( voir Théorème \ref{thm:324} et Théorème \ref{thm:326}) . Pour mémoire , la
base symplectique $( S_j  , \alpha_j )$ décrite précédemment est une base orthnormée
de $( $ G $ ,\omega , \hat{J})$ . En appliquant le théorème \ref{thm:401} on déduit
que la variété $\Gamma \backslash G$ est difféomorphe à un quotient fin du tore
plat. Soit $ \Gamma_o $ le sous-groupe abélien de $ \Gamma $ comme dans la preuve du
théorème \ref{thm:402}. Les variétés $ \Gamma \backslash G $ et $ \Gamma_
o\backslash G $ ont les mêmes nombres de BETTI. La courbure de Rici de $(  G ,
\omega , J )$ étant nulle un théorème de BOCHNER assure que la dimension de
l'algèbre des champs de Killing de $ (\Gamma \backslash G, g_{\hat{J} } ) $ est
égale au premier nombre de BETTI de la variété $\Gamma \backslash G $. En vertu du
théorème HATTORI ( ou de RAGHUNATHAN ) le premier nombre de BETTI de $\Gamma
\backslash G $ est égal à la dimension du premier espace de cohomologie réelle de
l'algèbre de Lie des champs de vecteurs invariants à gauche dans le groupe de Lie $
G $. Compte du théorème \ref{thm:402} de ce travail le premier nombre de BETTI de
$\Gamma \backslash G $ est égal à $ 2m = \dim G$. Pailleurs la dimension de premier
espace cohomologie de Koszul-de Rham de G est égale à la codimension du sous-groupe
des commutateurs de $G$. Cela permet de conclure que le groupe G est commutatif ; il
s'en suit que la variété $\Gamma \backslash G $ est difféomorphe au tore plat. \\

{\it Remarque}. Naturellement le théorème \ref{thm:411} est une généralisant non
triviale du théorème des tores Kählériens démontré par C. Benson et C Gordon d'une
part , [ BG1 ], ( voir aussi [McD1] )  . Nous observons que dans le cas des
nilvariétés traité par ces auteurs ,leur resultat peut être déduit directement de
notre théorème de réduction Kählérienne , Théorème \ref{thm:323}. En effet si on se
place dans les hypothèses du théorème \ref{thm:323} et si on suppose en outre que le
groupe de Lie $ G $ en jeu est nilpotent alors on peut choisir la configuration
(\ref{equation:8}) de sorte que le réseau $\Gamma$ donne lieu à des réseaux
cocompacts dans les sous-groupes de Lie $ G_1 $ et $ G_{2m+1} $ . Ainsi la variété
réduite $( G_{2m} , \omega )$ héritera à son tour d'un réseau cocompact
$\Gamma_{2m}$  , voir Remarque 2.2.3 [NB 6]. Par passage au quotient modulo les
réseaux le théorème de réduction Kählérien , Théorème \ref{thm:323}, determine la
variété Kählérienne compacte $( \Gamma_{2m} \backslash G_{2m} , \omega ,J_0 )$ .

La démonstration du théorème 1 de [ BD1 ] peut alors être faite par recurrence sur
la dimension des nilvariétés . Ainsi une fois supposé par hypothèse de recurrence
que $\Gamma_{2m} \backslash G_{2m} $ est un tore plat, le diagramme
(\ref{equation:8}) donne lieu à un diagramme analogue suivant :

$$
\begin{array}{cc}
{\mathbb T}^1  & \\
\downarrow & \\
M_{2m+1} & \longrightarrow \Gamma \backslash G \longrightarrow {\mathbb T}^1 \\
\downarrow & \\
{\mathbb T}^{2m}  &
\end{array}
$$

dans lequel la colonne et la ligne sont des fibrés en tores plats . En vertu du
Théorème \ref{thm:322} la colonne est un fibré trivial . Par conséquent la variété $
M_{2m+l} $ est difféomorphe au tore plat . Il reste seulement montrer que sous les
hypotèses du théorème 1 de [BGl] la ligne du diagramme ci-dessus est une fibration
triviale. Il suffit pour cela de combiner le théorème de Nomizu avec la suite
spectrale de Hochschild-Serre. \par

\noindent On peut exprimer le théorème 4.0.1 sous une forme qui fasse apparaître les
résultats de [BG1] et de [McD1] comme corollaires de ceux obtenus dans ce travail.

\begin{thm}
Soient $ \omega  \in \omega^ 2_l( G) $  une forme symplectique et $ \Gamma \subset G
$  un réseau cocompact dans le groupe de Lie complètement résoluble $ G\, . $  Si $
(\Gamma \backslash G,\omega) $ possède une structure k\"ahlérienne $ (\Gamma
\backslash G,\omega ,J) $ alors $ G $ est commutatif
\end{thm}

\noindent Preuve. On sait que les variétés compactes $ \Gamma_ o\backslash G $ et $
\Gamma \backslash G $ ont les mêmes nombres de BETTI réels.
\par \noindent En vertu d'un classique théorème de Hurewicz on a les égalités
$$ b_1(\Gamma \backslash G) = {\rm rang}({\mathcal D}\Gamma \backslash \Gamma)  =
{\rm rang} (\Gamma_ o ) = \dim G $$ o\`u $ {\mathcal D}\Gamma $ est le sous-groupe
des commutateurs de $ \Gamma \, . $
\par
\noindent Puisque $ b_1(\Gamma \backslash G) = \dim H^1(g,\bkR) $ o\`u $\frak g$ est
l'algèbre de Lie de $ G $ il en résulte que $ {\mathcal D}{\frak g} = [{\frak
g},{\frak g}] = (0).$
\par
\smallskip
\noindent Pour terminer on va donner un exemple de groupe de Lie résoluble qui
contient un réseau uniforme abélien sans être complètement résoluble et en tirer la
réponse à une conjecture de [BG 2].

\begin{example}
\label{exo:413}
Soit $ G $ le groupe de Lie dont la variété analytique sous-jacente
est $ {\mathbb C}^2 $ muni du produit

\begin{equation}
\label{equation:12} (z_1,z_2 ) (z^{\prime}_ 1,z^{\prime}_ 2 ) = (z_1+ z^{\prime}_ 1
\exp ( \sqrt{ -1} {\pi \over 2} y_2 ),z_2+z^{\prime}_ 2 )
\end{equation}

ou $ z_i= x_i+ \sqrt{ -1} y_i\quad i = 1,2\, . $ \par \noindent L'expression
(\ref{equation:12}) définit dans $ {\mathbb C}^2 $ une structure de groupe de Lie
réel de dimension 4. On pose $ \Gamma  = [{\mathbb Z} + \sqrt{ -1} {\mathbb Z} ]
\times [{\mathbb Z} + 4 \sqrt{ -1} {\mathbb Z} ] $ o\`u $ {\mathbb Z} $ est l'anneau
des nombres entiers ; $ \Gamma $ est un réseau commutatif cocompact dans $ G\, . $
Les translations à gauche par les $ \Gamma  \in \Gamma $ coïncident avec les
translations dans le groupe vectoriel $ \bkR^ 4, $ de sorte que la variété quotient
$ \Gamma \backslash G $ est le tore plat $ {\mathbb T}^ 4. $ Soit $ h $ la métrique
k\"ahlérienne canonique de $ {\mathbb C}^ 4\, : $

$$ h(z) = dz_1d\bar  z_1+ dz_2d\bar  z_2\, . $$

La métrique $ h $ est invariante par les translations à gauche par les $ \Gamma \in
\Gamma \, ; $ elle donne par passage au quotient la métrique de Kähler standard du
tore plat $ {\mathbb T}^ 4. $

\par \noindent Pour $ z \in  G $ $ Ad(z) $ est l'image de $ z $
par la représentation adjointe. Dans la base canonique de $ {\mathbb C}^ 2 $ la
matrice de $ Ad(z) $ est

$$ Ad(z) =  \begin{pmatrix}  \cos {\pi y_2 \over 2} &
  - \sin {\pi y_2 \over 2} &   0 &   - {\pi y_1
\over 2} \\  \sin {\pi y_2 \over 2} &   \cos {\pi y_2 \over 2} &   0 &
 { \pi x_1 \over 2} \\  0 &   0 &
1 &   0 \\  0 &   0 &   0 &   1 \end{pmatrix}
$$

les valeurs propres des $ Ad(z) $ ne sont pas toutes réelles. Ainsi $ G $ n'est pas
complètement résoluble. \par \noindent Il a été conjecturé ce qui suit, [BG 2].

\noindent Soit $ G $ un groupe de Lie résoluble contenant un réseau cocompact $
\Gamma $ tel que $ \Gamma \backslash G $ est difféomorphe au tore plat. \par
\noindent Si $ \Gamma \backslash G $ possède une métrique k\"ahlérienne alors $ G $
est complètement résoluble.
\end{example}

L'exemple \ref{exo:413}  montre qu'il n'en est rien.

\subsection{Groupes de Lie Kahhlériens}
 Nous rappelons un résultat déjà ancien [ NB7 ] qui peut être déduit
sans difficultés du théorème de réduction Kählérien . Un groupe de Lie muni d'une
structure Kählérienne invariante par les translations à gauche est appelé groupe de
Lie Kählérien . Le cardre naturel de traitement des groupes de Lie Kählériens est
celui des variétés Kählériennes homogènes , [ DV ]  [ DN ] . Cependant les groupes
de Lie Kählériens complètement résolubles jouissent des propriétés assez fines qui
viennent d'un avatar de la géométrie Riemannienne , [ DP ] qu'est la géométrie
hessienne . Soit $( M , D )$ une variété localement plate .Une métrique Riemannienne
$g$ définie dans M est dite localement hessienne si au voisinage de tout point de
$M$ existe une fonction lisse $h$ dont la $D$-hessienne est le tenseur métrique $g$,
c'est à dire que l'on a $g = D^2(h)$ . Les groupes de Lie Kählériens sont en fait
des variété tube des groupes de Lie héssiens . \\

{\it Modèles standard}. Soit $G$ un groupe de Lie connexe , simplement connexe dont
l'algèbre de Lie des champs des vecteurs invariants par les translations à gauche
est $\frak g $ . Soit ${\frak g}^{\ast} $ l'espace vectoriel dual de $\frak g $ ; si
$(G,D)$ est une structure localement plate invariante à gauche dans $G$ alors
l'espace vectoriel ${\frak g} + {\frak g}^{\ast} $ en hérite d'une structure
d'algèbre de Lie dont le crochet est défini par la formule $   [ ( f , X ) , ( f'
,X' ) ]  =( fD_X' - f'D_X , [ X , X' ] ) $. Le fibré cotangent T*G  possède une
structure de groupe de Lie dont l'algèbre de Lie est identifiée avec ${\frak
g}^{\ast} + {\frak g} $ , nous désignerons ce groupe de Lie par $G_D$ . Soit
$\omega_D$ la 2-forme différentielle invariante sur $G_D$ engendrée par
 $\omega (( f,X ) , ( f' , X )) = f( X' ) - f' ( X ) $.
Cette forme est fermée . D'un autre côté $G_D$ possède une structure de variété
localement plate définie par la connection invariante à gauche suivante $\hat{D}_{(
f ,X )} ( f' , X' ) = ( - f' D_X  , D_X X'  )$ . Le triplet $( T^{\ast}G
,\omega,\hat{D} )$ définit une structure de variété symplectique localement plate
invariante par les translations à gauche dans $T^{\ast}G$ . Supposons maintenant que
$G$ possède une métrique Riemannienne D-hessienne , $\frak g $ ,alors $T^{\ast}G$
possède une métrique de Kähler $h$ dont la partie réelle est quasi-fibré , [ DP ] .
Notons $f$ l'isomorphisme de $\frak g$ sur ${\frak g}*$ défini par la metrique $g$.
Alors $T^{\ast}G$ possède une structure de variété analytique complexe dont le
tenseur
 de la structure presque complexe $J$ est défini par $J ( f(Y) , X ) = f( X ) , - Y )$
. On vérifie par un calcul direct que la forme symplectique $ \omega $ est
J-invariante est que $ \omega ( JX ,X )$ est définie positive. Si on suppose que la
métrique $D$-hessienne $ g $ est invariante à gauche dans G alors le triplet
$(T^{\ast}G ,\omega, J)$ définit un groupe de Lie Kählérien dont la métrique
Riemannienne est quasi-fibré sur $(G,g)$ . A ce point on n'a pas fait d'hypothèse
particulière sur $G$ ; si ce dernier est complètement résoluble alors le modèle qui
vient d'être construit est standard. Plus précisément on a

\begin{thm}Tout groupe de Lie Kählérien connexe simplement connexe et complètement
résoluble est le cotangent d'un  groupe de Lie complètement résoluble $G$ muni d'une
structure structure affine invariante à gauche $( G , D )$ et d'une métrique
Riemannienne localement $D$-hessienne invariante à gauche .
\end{thm}

{\it Esquisse de démonstration}. Soit  $( G ,  \omega  , J )$ un groupe de Lie
Kählérien de dimension réelle $2m$ . Grâce au théorème 3.2.2 de \cite{NB 6} on sait
que $ ( G ,\omega )$ possède un feuilletage lagrangien bi-invariant $L \subset TG$.
Le théorème \ref{thm:324} assure que l'orthogonal Riemannien de $ L , J( L ) $ ,
définit un système différentiel complètement intégrable. Notons $G_L$ et $G_{J(L)}$
respectivement les feuilles de $ L $ et de $ J (L ) $ qui passent par l'élément
neutre de $ G $; $ G $ est produit semi-direct de $ G_{J(L)} $ par $ G_L $. Le
tenseur de la structure presque complexe J définit un isomorphisme canonique de
l'agèbre de Lie de $ G_L $ sur son espace vectoriel dual ; on peut ainsi identifié
$G$ avec le fibré cotangent $ T^{\ast}G_{J(L)} $. On définit dans $ G_{J(L)} $ la
connexion linéaire invariante à gauche dont la dérivation covariante $D$ est donnée
par la formule $ D_X Y = - J [ X , J(Y) ] $ , $X$ et $Y$ étant des éléments de
l'algèbre de Lie de $ G_{J(L)} $ .Cette connexion est localement plate . Sous
l'identification de $G$ avec $T^{\ast}G_{J(L)}$ la forme symplectique $\omega $
dévient $ \omega_D $. Pour terminer on vérifie directement la véracité de
l'assertion suivante.

\begin{prop}
La restriction à $ G_{J(L)} $ de la partie réelle de h est localement D-hessienne.
\end{prop}

Pour terminer on observe que pour un groupe de Lie complètement résoluble non
commutatif G les deux  assertions suivantes s'excluent mutuellement : ( a )  $G$
possède une structure de groupe de Lie Kählérien , ( b )  $G$ contient un réseau
cocompact .

\bibliographystyle{amsplain}

\begin{thebibliography}{ABCDEF}


\bibitem[ABCKT]{ABCKT}
J. AMOROS, N. BURGER, K. CORLETTE, D. KOSCHICK and D. TOLEDO. {\it Fundamental
groups of compact Kähler manifolds.} Amer Math. Soc. Math. survey and Monographs
vol. 44.

\bibitem[AU 1]{AU 1}
L. AUSLANDER. {\it Bieberbach's groups on space groups and discrete uniform
subgroups of Lie groups}. Ann. Math. 71 (1960) 579--590.



\bibitem[AU 2]{AU 2}
L. AUSLANDER. {\it An exposition of the structures of solvmanifolds.} Bull. of Amer.
Math. Soc. 79 (1973) 227--261.

\bibitem[BG 1]{BG 1}
C. BENSON and C. GORDON. {\it Kähler and symplectic structures on nilmanifolds}.
Topology 27 (1988) 513--518.

\bibitem[BG 2]{BG 2}
C. BENSON and C. GORDON. {\it Kähler structures on compact solvmanifolds.} Proc. of
the Amer. Math. Soc. 108 (1990) 971--990.

\bibitem[BE]{BE}
A. BESSE. Enstein manifolds. Ergebnisse der Math. und ihrer Grenzgebiete 3. Folge.
Band 10.

\bibitem[B0]{B0}
N. BOURBAKI.  Groupes et Algèbres de Lie. vol. 7. Hermann Paris.

\bibitem[BW]{BW}
F. BRUHAT and H. WHITNEY. Comment. Math. Helv. 33 (1959) 132--160.

\bibitem[CFG]{CFG}
L. CORDERO, M. FERNANDEZ and A. GRAY. {\it Symplelctic manifolds with no Kähler
structure.} Topology 25 (1986) 375--380.

\bibitem[DGMS]{DGMS}
P. DELIGNE, P. GRIFFITHS, J. MORGAN and D. SULLIVAN. {\it Real homotopy theory of
Kähler manifolds.} Inv. Math. 29 (1975) 245--274.

\bibitem[DN]{DN}
J. DORFNEISTER and K. NAKAJIMA. {\it The fundamental conjecture.} Acta Math. 160
(1988) 13--70.

\bibitem[DP]{DP}
Ph. DELANOE , {\it Remarques sur les variétés localement hessiennes}. Osaka J. Math.
26 (1989), 65-69.

\bibitem[FGG]{FGG}
M.FERNANDEZ , M. GOTAY and A. GRAY {\it Compact parallelizable four dimensional
symplecric manifolds.} Proc. Amer Math Soc 103 ( 1988) 1209- 1212 .

\bibitem[FGH]{FGH}
D. FRIED, W. GOLDMAN and M.W. HIRSCH. {\it Affine manifolds with nilpotent holonomy
group.} Comment Math. Helv. 56 (1981) 487--523.

\bibitem[G]{G}
M. GROMOV. {\it  Sur le groupe fondamental d'une variété Kählérienne.} CR. Acad.
Sci. Paris 308 (1989) 67--70.

\bibitem[GPV 1]{GPV 1}S.G. GINDIKIN, I.I.
PIATECCKII, SAPIRO and E.B. VINBERG. {\it Homogeneous K\"ahler manifolds.} Geometry
of homogeneous bounded domains, CIME 1967.

\bibitem[GPV 2]{GPV 2} S.G. GINDIKIN, I.I.
PIATECCKII, SAPIRO and E.B. VINBERG. {\it On the classification and canonical
realization of complex homogeneous bounded domains.} Trans. Moscow Math. Soc. 13
(1964) 340--403.

\bibitem[GH]{GH}
W. GOLDMAN and M.W. HIRSCH. {\it Affine manifolds and Orbits of Algebraic groups.}
Trans. Amer. Math. Soc. 295 (1986) 175--197.

\bibitem[GM]{GM}
W. GOLDMAN and J.J. NILLSON. {\it The deformation theory of representations of
fundamental groups of compact Kähler manifolds.} Publ. Math. IHES 67 (1988) 43--96.

\bibitem[GS]{GS}
V. GUILLEMIN and M. STENZEL. {\it  Grauert tube and the homogeneous Monge-Ampere
equations.} J. Diff. Geom. 34 (1991) 561--570.

\bibitem[GSS]{GSS}
V.GULLEMIN and S. STERNBERG. {\it The algebraic model for transitive differential
geometry.} Bull. Amer Math Soc 70 ( 1964 ) 16-47 .

\bibitem[HA]{HA}
A. HATTORI. {\it Spectral sequence in the de RHAM cohomology of fibre bundles.} Fac.
sei. Univ. Tokyo 8 (1960) 289--331.

\bibitem[HH]{HH}
A.HESS. {\it Connection on symplectic manifolds and Geometric quantization.} Lecture
Notes in Math.

\bibitem[HW]{HW}
W.Y. HISIANG. {\it Cohomology theory of topological transformation groups}. Erg. der
Math. une chrer Gr. band 85 Springer--Verlag.

\bibitem[JD]{JD}
L.D. JOHNSON. {\it Kähler submersion and holomorphic connection.} Jour Diff Geom  15
(1980) 71-79.

\bibitem[JW]{JW}
L.D.JOHNSON  and L.B. WITT . Totally geodesic foliations . Jour Diff Geom 15
(1980 ) 225-235 .


\bibitem[KL 1]{KL 1}
J.L. KOSZUL. {\it Homologie et Cohomologir des Algèbres de Lie.} Bull. Soc. Math.
France 48 (1950) 65--127.

\bibitem[KL 2]{KL 2}
J.L. KOSZUL. {\it Déformations des variétés localement plates.} Ann. Inst. Fourier
18 (1968).

\bibitem[L]{L}P. LIBERMANN. {\it  Problème
d'équivalence en géométrie symplectique.} Astérisque 107 (1983) 43--68. \par

\bibitem[MA]{MA}
Y. MATSUSHIMA. {\it  On discrete subgroups and Homogeneous spaces of nilpotent Lie
groups.} Nagoya Math. Jour. 2 (1951) 95--110.

\bibitem[McD 1]{Mc D 1}
Dusa Mc DUFF. {\it The moment map for circle action on symplectic manifolds} JGP 5
(1988) 149--160.

\bibitem[McD 2]{Mc D 2}
Dusa Mc DUFF. {\it Example of simply connected symplectic manifolds non
k\"ahlerian.} Jour. Diff. Geom. 70 (1984) 267--277.

\bibitem[MI]{MI}
J. MILNOR. {\it Fundamental groups of complete affinely flat manifolds.} Adv. in
Math. 25 (1977) 178--187.

\bibitem[MO 1]{MO 1}
G.D. MOSTON. {\it Cohomology of topological groups and solomanifolds.} Ann. of Math.
73 (1961) 20--48.

\bibitem[MO 2]{MO 2}
G.D. MOSTON. {\it  Fundamental groups of homogeneous spaces.} Ann. of Math. 66
(1957) 249--255.

\bibitem[MW]{MW}
J. MARSDEN and WEINSTEIN. {\it Reduction of symplectic manifolds with symmetry.}
Reports on Math. Phys. 5 (1974) 121--130.

\bibitem[NB 1]{NB 1}
NGUIFFO B. BOYOM. {\it Variétés symplectiques affines.} Manuscripta Math. 64 (1989)
1--33.

\bibitem[NB 2]{NB 2}
NGUIFFO B. BOYOM. {\it Structures affines homotpes à zéro.} Jour. Diff. Geom. 33
(1990) 859--911.

\bibitem[NB 3]{NB 3}
NGUIFFO B. BOYOM. {\it Structures affines isotropes des groupes de Lie.} Annalli
Scuo. Norm. Sup. Pisa.

\bibitem[NB 4]{NB 4}
NGUIFFO B. BOYOM. {\it Métriques k\"ahlériennes affinement plates.} Proc. London
Math. Soc. 3 (1993) 358--380.

\bibitem[NB 5]{NB 5}
NGUIFFO B. BOYOM. {\it Structures localement plates dans certaines variétés
symplectiques.} Math Scand. 76 (1995) 61--84.

\bibitem[NB 6]{NB 6}
NGUIFFO B. BOYOM. {\it Sur les actions hamiltoniennes quasi-primitives} ( soumis ).


\bibitem[NK]{NK}
K. NOMIZU. {\it  On the cohomology of compact homogeneous spaces of nilpotent Lie
groups.} Ann. of Math. 59 (1954) 531--538.

\bibitem[ON]{ON}A.L. ONISCHIK. {\it Inclusion
relations among transitive compact groups.} Amer. Math. Soc. Translations S.2 50
(1966) 5--58.

\bibitem[RA]{RA}S. RAGHUNATHAN.  Discrete
subgroups of Lie groups. Ergebnisse der Math. 68 Springer Verlag 1972.

\bibitem[RC]{RC} C. RAMANUJAM.{\it  A topological
characterisation of the affine plane as an Algebraic Variety.} Ann. of Math. 94
(1971) 69--88.

\bibitem[RW]{RW}J. ROELS and A. WEINSTEIN. {\it
On functions whose Poisson brackets are constant.} J. Math. Phys. 12. 1482--1486.

\bibitem[SA]{SA}M. SATO. {\it Sur certains groupes
de Lie résolubles.} Sci. pap. Univ. Tokyo 7 (1957) $I$--$II$ ; 157--168.

\bibitem[SE]{SE}J.P. SERRE. GAGA. Ann. Inst.
Fourier 6 (1956) 1--42.

\bibitem[SH]{SH}
H. SHIMA. {\it Homogeneous hessian manifolds}. Ann. Inst. fourier 30 (1980), 91-128.

\bibitem[SS]{SS} I.M. SINGER and S. STERNBERG. {\it The infinte groups of Lie and Cartan.}
Journ. d'Ananlyse Math. Jerusalem 15 ( 1965 ) 1-114 .

\bibitem[TH]{TH}W. THURSTON. {\it Some simple of
symplectic manifolds.} Proc. Amer. Math. Soc. 55 (1976) 456--478.

\bibitem[WE]{WL}A. WEILL.  Introduction à
l'étude des variétés k\"ahlériennes. Hermann Paris 1958.

\bibitem[WN]{WN}A. WEINSTEIN.{\it Symplectic
manifolds and their lagrangian submanifolds.} Adv. In Math. 6 (1971) 329--346.

\bibitem[WOJ]{WOJ} J.A. WOLF . Spaces of constant curvature , McGraw-Hill , New York , 1967 .

\bibitem[WO]{WO}N.M.J. WOODHOUSE. Geometric
quantization. Oxford Math. Monographs. Larendon Press Oxford.

\end{thebibliography}

\end{document}